\documentclass[10pt,twocolumn,twosde]{IEEEtran}
\usepackage{amsmath, amssymb, a4wide}
\usepackage{balance}
\usepackage{amsthm}
\usepackage{graphicx}
\usepackage{bm}
\usepackage{relsize}
\usepackage{booktabs}
\usepackage{dsfont}
\usepackage{url}
\usepackage{enumerate}
\usepackage[latin1]{inputenc}
\usepackage{multirow}
\usepackage[numbers,sort&compress]{natbib}
\usepackage{caption}
\usepackage{subcaption}
\usepackage{array}
\newcolumntype{L}[1]{>{\raggedright\let\newline\\\arraybackslash\hspace{0pt}}m{#1}}
\newcolumntype{C}[1]{>{\centering\let\newline\\\arraybackslash\hspace{0pt}}m{#1}}
\newcolumntype{R}[1]{>{\raggedleft\let\newline\\\arraybackslash\hspace{0pt}}m{#1}}
\usepackage{enumerate}
\usepackage{color}
\usepackage[colorlinks,citecolor=blue,urlcolor=blue]{hyperref}
\usepackage[T1]{fontenc}
\usepackage{setspace}
\usepackage{mathrsfs}
\usepackage[margin=0.8in]{geometry}
\usepackage{hyperref}
\usepackage{textcomp}

\DeclareMathOperator\tr{Tr}

\newtheoremstyle{definition}%
  {0.5em}
  {\topsep}
  {\upshape}
  {}
  {\bfseries}
  {.}
  { }
  {}

\theoremstyle{plain}
\newtheorem{theorem}{Theorem}
\newtheorem{proposition}{Proposition}
\newtheorem{lemma}{Lemma}
\newtheorem{corollary}{Corollary}
\newtheorem{remark}{Remark}
\newtheorem{definition}{Definition}
\theoremstyle{definition}
\newtheorem{assumption}{Assumption}

\allowdisplaybreaks[2]

\begin{document}
\title{Large-dimensional behavior of regularized Maronna's M-estimators of covariance matrices}

\author{N. Auguin$^{\star}$, D. Morales-Jimenez$^{\dagger}$, M. R. McKay$^{\star}$, R. Couillet$^{\ddagger}$

\thanks{$^{\star}$N. Auguin and M. R. McKay are with the Department of Electronic and Computer Engineering, Hong Kong University of Science and Technology, Clear Water Bay, Kowloon, Hong Kong. E-mail: nicolas.auguin@connect.ust.hk, m.mckay@ust.hk.}
\thanks{$^{\dagger}$D. Morales-Jimenez is with the Institute of Electronics, Communications and Information Technology, Queen's University Belfast, Belfast BT3 9DT, United Kingdom. E-mail: d.morales@qub.ac.uk.}
\thanks{$^\ddagger$R. Couillet is with GIPSA-lab at Universit\'e Grenoble-Alpes, St-Martin d'H\`eres, France, and with L2S, CentraleSup\'elec at Universit\'e Paris-Saclay, Gif-sur-Yvette, France. E-mail:  romain.couillet@gipsa-lab.grenoble-inp.fr, romain.couillet@centralesupelec.fr.}
\thanks{N. Auguin, D. Morales-Jimenez and M. R. McKay were supported by the Hong Kong RGC General Research Fund under grant numbers 16206914 and 16203315. R. Couillet's work was supported by the ANR Project RMT4GRAPH (ANR-14-CE28-0006). }}

\maketitle

\begin{abstract}
Robust estimators of large covariance matrices are considered, comprising regularized (linear shrinkage) modifications of Maronna's classical M-estimators. These estimators provide robustness to outliers, while simultaneously being well-defined when the number of samples does not exceed the number of variables. By applying tools from random matrix theory, we characterize the asymptotic performance of such estimators when the numbers of samples and variables grow large together. In particular, our results show that, when outliers are absent, many estimators of the regularized-Maronna type share the same asymptotic performance, and for these estimators we present a data-driven method for choosing the asymptotically optimal regularization parameter with respect to a quadratic loss. Robustness in the presence of outliers is then studied: in the non-regularized case, a large-dimensional robustness metric is proposed, and explicitly computed for two particular types of estimators, exhibiting interesting differences depending on the underlying contamination model. The impact of outliers in regularized estimators is then studied, with interesting differences with respect to the non-regularized case, leading to new practical insights on the choice of particular estimators. 
\end{abstract}


\section{Introduction}
\label{sec:intro}
Covariance or scatter matrix estimation is a fundamental problem in statistical signal processing \cite{abramovich2007diagonally,pascal2014generalized}, with applications ranging from wireless communications \cite{tulino2004random} to financial engineering \cite{ledoit2003improved} and biology \cite{schafer2005shrinkage}. Historically, the sample covariance matrix (SCM) $\frac{1}{n}\sum_{i=1}^{n}\mathbf{y}_i\mathbf{y}_i^{\dagger}$, where $\mathbf{y}_1, \cdots, \mathbf{y}_n \in \mathbb{C}^N$ are zero-mean data samples, has been a particularly appealing choice among possible estimators. The SCM is known to be the maximum likelihood estimator (MLE) of the covariance matrix when the $\mathbf{y}_i$ are independent, identically distributed zero-mean Gaussian observations, and its simple structure makes it easy to implement. Nonetheless, the SCM is known to suffer from three major drawbacks: first, it is not resilient to outliers nor samples of impulsive nature; second, it is a poor estimate of the true covariance matrix whenever the number of samples $n$ and the number of variables $N$ are of similar order; lastly, it is not invertible for $n<N$. The sensitivity to outliers is a particularly important issue in radar-related applications \cite{ward1981compound,billingsley1999statistical}, where the background noise usually follows a heavy-tailed distribution, often modeled as a complex elliptical distribution \cite{kelker1970distribution,ollila2012complex}. In such cases, the MLE of the covariance matrix is no longer the SCM. On the other hand, data scarcity is a relevant issue in an ever-growing number of signal processing applications where $n$ and $N$ are generally of similar order, possibly with $n<N$ \cite{ledoit2003improved,schafer2005shrinkage,mestre2008modified,nadler2010nonparametric}.
New improved covariance estimators are needed to account for both potential data anomalies and high-dimensional scenarios.

In order to harness the effect of outliers and thus provide a better inference of the true covariance matrix, robust estimators known as M-estimators have been designed \cite{huber1964robust,maronna1976robust,tyler1987distribution,ollila2012complex}. Their structure is non-trivial, involving matrix fixed-point equations, and their analysis challenging. Nonetheless, significant progress towards understanding these estimators has been made in large-dimensional settings \cite{chen2011robust,couillet2013robust,couillet2014large,zhang2014marchenko,morales2015large}, motivated by the increasing number of applications where $N,n$ are both large and comparable. Salient messages of these works are: (i) outliers or impulsive data can be handled by these estimators, if appropriately designed (the choice of the specific form of the estimator is important to handle different types of outliers) \cite{morales2015large}; (ii) in the absence of outliers, robust M-estimators essentially behave as the SCM and, therefore, are still subject to the data scarcity issue \cite{couillet2013robust}.

To alleviate the issue of scarce data, regularized versions of the SCM have originally been proposed \cite{abramovich1981controlled,carlson1988covariance}. Such estimators consist of a linear combination of the SCM and a shrinkage target (often the identity matrix), which guarantees their invertibility, and often provides a large improvement in accuracy over the SCM when $N$ and $n$ are of the same order. Nevertheless, the regularized SCM (RSCM) inherits the sensitivity of the SCM to outliers/heavy-tailed data. To alleviate both the data scarcity issue and the sensitivity to data anomalies, regularized M-estimators have been proposed
 \cite{abramovich2007diagonally,chen2011robust,pascal2014generalized,ollila2014regularized}.  Such estimators are similar in spirit to the RSCM, in that they consist of a combination of a robust M-estimator and a shrinkage target. However, unlike the RSCM, but similar to the estimators studied in \cite{couillet2013robust,morales2015large}, these estimators are only defined implicitly as the solution to a matrix fixed-point equation, which makes their analysis particularly challenging. 

In this article, we propose to study these robust regularized estimators in the double-asymptotic regime where $N$ and $n$ grow large together. Building upon recent works \cite{couillet2014large,couillet2015random}, we will make use of random matrix theory tools to understand the yet-unknown asymptotic behavior of these estimators and subsequently to establish design principles aimed at choosing appropriate estimators in different scenarios. In order to do so, we will first study the behavior of these regularized M-estimators in an outlier-free scenario. In this setting, we will show that, upon optimally choosing the regularization parameter, most M-estimators perform asymptotically the same, meaning that the form of the underlying M-estimator does not impact the performance of its regularized counterpart in clean data. Second, we will investigate the effect of the introduction of outliers in the data, under different contamination models. Initial insights were obtained in \cite{morales2015large} for non-regularized estimators, focusing on the weights given by the M-estimator to outlying and legitimate data. However, the current study, by proposing an intuitive measure of robustness, takes a more formal approach to qualify the robustness of these estimators. In particular, we will demonstrate which form of M-estimators is preferable given a certain contamination model, first in the non-regularized setting, and then for regularized estimators. \newline  

\emph{Notation:} $||\mathbf{A}||$, $||\mathbf{A}||_F$ and $\tr \mathbf{A}$ denote the spectral norm, the Frobenius norm and the trace of the matrix $\mathbf{A}$, respectively. The superscript $(\cdot)^\dagger$ stands for Hermitian transpose. Thereafter, we will use $\lambda_1(\mathbf{A}) \leq \cdots \leq \lambda_N(\mathbf{A})$ to denote the ordered eigenvalues of the square matrix $\mathbf{A}$. The statement $\mathbf{A} \succ 0$ (resp.\@ $\succeq 0$) means that the symmetric matrix $\mathbf{A}$ is positive definite (resp.\@ positive semi-definite). The arrow $\xrightarrow{\substack{\mathrm{a.s.}}}$ designates almost sure convergence, while $\mathbf{\delta}_x$ denotes the Dirac measure at point $x$. 

\section{Review of the large dimensional behavior of non-regularized M-estimators}
\label{sec:2}
\subsection{General form of non-regularized M-estimators}
In the non-regularized case, robust M-estimators of covariance matrices are defined as the solution (when it exists) to the equation in $\mathbf{Z}$ \cite{maronna1976robust}
\begin{align}
\label{eq:implicit Z}
\mathbf{Z} = \frac{1}{n} \sum_{i=1}^n u\left(\frac{1}{N}\mathbf{y}_i^\dagger \mathbf{Z}^{-1} \mathbf{y}_i \right) \mathbf{y}_i \mathbf{y}_i^\dagger,
\end{align}
where $\mathbf{Y}=[\mathbf{y}_1, \cdots, \mathbf{y}_{n}] \in \mathbb{C}^{N \times n}$ represents the data matrix, and where $u$ satisfies the following properties:
\begin{itemize}
\item $u$ is a nonnegative, nonincreasing, bounded, and continuous function on $\mathbb{R}^+$, 
\item $\phi : x \mapsto x u(x)$ is increasing and bounded, with $\phi_\infty \triangleq \lim_{x \rightarrow \infty} \phi(x) > 1$.
\end{itemize}
If well-defined, the solution of (\ref{eq:implicit Z}) can be obtained via an iterative procedure (see, for example, \cite{kent1991redescending,couillet2013robust}). Intuitively, the $i$-th data sample is given a weight $u(\frac{1}{N}\mathbf{y}_i^\dagger \mathbf{Z}^{-1} \mathbf{y}_i )$, which should be smaller for outlying samples than for legitimate ones. The choice of the $u$ function determines the degree of robustness of the M-estimator. As a rule of thumb, the larger $\phi_\infty$, the more robust the underlying M-estimator to potential extreme outliers \cite{maronna1976robust}. However, such increased robustness is usually achieved at the expense of accuracy \cite{ollila2012complex}.

A related and commonly-used estimator is Tyler's estimator \cite{tyler1987distribution}, which is associated with the unbounded function {\color{black}$u_\mathrm{Tyler}(x) = 1/x$}. {\color{black} Recent papers dealing with the application of Tyler's estimator in engineering include \cite{soloveychik2014tyler,yang2015robust,couillet2016second}.} We remark however that, for such $u$ function, the existence of a solution to (\ref{eq:implicit Z}) depends on the data at hand (see, e.g., \cite{tyler1987distribution,sun2014regularized,ollila2014regularized}). To avoid this issue, {\color{black}we here focus on a class of estimators with bounded $u$ functions} (as prescribed above). Examples of practical interest, which we study in some detail, include
\begin{align}
\label{u_tyler}
u_{\mathrm{M-Tyler}}(x) &\triangleq K \frac{1+t}{t+x} \\
\label{u_huber}
u_{\mathrm{M-Huber}}(x) &\triangleq K \min\left\{1,\frac{1+t}{t+x}\right\},
\end{align}
for some $t,K>0$. For a specific $t$, $u_{\mathrm{M-Tyler}}$ is known to be the MLE of the true covariance matrix when the $\mathbf{y}_i$ are independent, zero-mean, multivariate Student vectors \cite{morales2015large}, whereas $u_{\mathrm{M-Huber}}$ refers to a modified form of the so-called Huber estimator \cite{huber2011robust}.  Observe that for these functions, $\phi_\infty = K(1+t)$, such that the robustness of the associated M-estimator to extreme outliers is controlled by both $t$ and the scale factor $K$. In what follows, with a slight abuse of terminology, we will refer to these estimators as ``Tyler's'' and ``Huber's'' estimators, respectively. 

\subsection{Asymptotic equivalent form under outlier-free data model}
Assume now the following ``outlier free'' data model: let $\mathbf{y}_i$ be $N$-dimensional data samples, drawn from $\mathbf{y}_i=\mathbf{C}_N^{1/2}\mathbf{x}_i$, where $\mathbf{C}_N \in \mathbb{C}^{N \times N} \succ 0$ is deterministic and $\mathbf{x}_1, \ldots, \mathbf{x}_{n}$ are random vectors, the entries of which are independent with zero mean, unit variance and finite $(8+\sigma)$-th order moment (for some $\sigma > 0$). With this model, we now recall the main result from \cite{couillet2013robust}.

\begin{theorem} \cite{couillet2013robust}
\label{th:Romain_no_shrinkage}
Assume that $c_N \triangleq N/n \rightarrow c \in (0, 1)$ as $N,n \rightarrow \infty$. Further assume that $ 0 < \liminf_N \{\lambda_1\left( \mathbf{C}_N \right) \} \leq \limsup_N \{\lambda_N\left( \mathbf{C}_N \right) \} < \infty$. Then, denoting by $\hat{\mathbf{C}}_N$ a solution to (\ref{eq:implicit Z}), we have
\begin{align*}
\left\lVert \hat{\mathbf{C}}_N - \hat{\mathbf{S}}_N \right\rVert\xrightarrow{\substack{\mathrm{a.s.}}} 0,
\end{align*}
where $\hat{\mathbf{S}}_N \triangleq \frac{1}{\phi^{-1}(1)} \frac{1}{n} \sum_{i=1}^n \mathbf{y}_i \mathbf{y}_i^\dagger$.
\end{theorem}

This shows that, up to a multiplying constant, regardless of the choice of $u$, Maronna's M-estimators behave (asymptotically) like the SCM. As such, in the absence of outliers, no information is lost. 

However, Theorem \ref{th:Romain_no_shrinkage} excludes the ``under-sampled'' case $N \geq n$. Regularized versions of Maronna's M-estimators have been proposed to alleviate this issue, in most cases considering regularized versions of Tyler's estimator ($u(x) =1/x$) \cite{abramovich2007diagonally,chen2011robust,sun2014regularized,pascal2014generalized}, the behavior of which has been studied in \cite{zhang2014marchenko,couillet2014large}. Recently, a regularized M-estimator which accounts for a wider class of $u$ functions has been introduced in \cite{ollila2014regularized}, but its large-dimensional behavior remains unknown. We address this in the next section. Moreover, note that Theorem~\ref{th:Romain_no_shrinkage} does not tell us anything about the behavior of different estimators, associated with different $u$ functions, in the presence of outlying or contaminating data. While progress to better understand the effect of outliers was recently made in \cite{morales2015large}, their study focused on non-regularized estimators. In this work, a new measure to characterize the robustness of different M-estimators will be proposed, allowing us to study both non-regularized and regularized estimators (Sections~\ref{sec:if} and \ref{sec:if_reg}).

\section{Regularized M-estimators: Large dimensional analysis and calibration}
\label{sec:shrinkage}
\subsection{General form of regularized M-estimators}
We consider the class of regularized M-estimators introduced in \cite{ollila2014regularized}, and given as the unique solution to 
\begin{align}
\begin{array}{ll}
\mathbf{Z} &\hspace{-2mm}=(1-\rho) \frac{1}{n}\sum_{i=1}^n u\left(\frac{1}{N}\mathbf{y}_i^\dagger\mathbf{Z}^{-1}\mathbf{y}_i\right)\mathbf{y}_i\mathbf{y}_i^\dagger+\rho \mathbf{I}_N
\end{array}, \label{implicit2}
\end{align} 
where $\rho \in (0,1]$ is a regularization (or shrinkage) parameter, and where $\mathbf{I}_N$ denotes the identity matrix. The introduction of a regularization parameter allows for a solution to exist when $N>n$. The structure of (\ref{implicit2}) strongly resembles that of the RSCM, defined as
\begin{align}
\label{LW}
\mathbf{R}(\beta) \triangleq (1-\beta)\frac{1}{n}\sum_{i=1}^{n}\mathbf{y}_i\mathbf{y}_i^{\dagger}+\beta \mathbf{I}_N,
\end{align}
where $\beta \in [0,1]$, also referred to as linear shrinkage estimator \cite{ledoit2004well}, linear combination estimator \cite{du2010fully}, diagonal loading \cite{carlson1988covariance}, or ridge regression \cite{hoerl1970ridge}. Regularized M-estimators are robust versions of the RSCM.

\vspace{-3mm}

\subsection{Asymptotic equivalent form under outlier-free data model}
Based on recent random matrix theory results, we now characterize the asymptotic behavior of these M-estimators. Under the same data model as that of Section \ref{sec:2}, we answer the basic question of whether (and to what extent) different regularized estimators, associated with different $u$ functions, are asymptotically equivalent. We need the following assumption on the growth regime and the underlying covariance matrix $\mathbf{C}_N$:
\begin{assumption}  
\hfill \begin{enumerate}[a.] 
\item $c_N \triangleq N/n \rightarrow c \in (0, \infty)$ as $N,n \rightarrow \infty$. 
\item $ \limsup_N \{\lambda_N\left( \mathbf{C}_N \right) \} < \infty$.
\item $\nu_n \triangleq \frac{1}{N}\sum_{i=1}^N \delta_{\lambda_i(\mathbf{C}_N)} $ satisfies $\nu_n \rightarrow \nu $ weakly with $\nu \neq \delta_0 $ almost everywhere.
\end{enumerate}
\end{assumption}
Assumption 1 slightly differs from the assumptions of Theorem \ref{th:Romain_no_shrinkage}. In particular, the introduction of a regularization parameter now allows $c\geq1$. Likewise, $\mathbf{C}_N$ is now only required to be positive semidefinite. 

For each $\rho \in (0,1]$, we denote by $\hat{\mathbf{C}}_N(\rho)$ the unique solution to (\ref{implicit2}). We first characterize its behavior in the large $n, N$ regime. To this end, we need the following assumption:
\begin{assumption}
$\phi_\infty = \lim_{x \rightarrow \infty} \phi(x) < \frac{1}{c}$. 
\end{assumption}
We now introduce an additional function, which will be useful in characterizing a matrix equivalent to $\hat{\mathbf{C}}_N(\rho)$. \\

\noindent \textbf{Definition.} Let Assumption 2 hold. Define {\color{black} $v:[0,\infty) \rightarrow (0,u(0)]$} as $v(x)=u(g^{-1}(x))$ where $g^{-1}$ denotes the inverse function of $g(x)=\frac{x}{1-(1-\rho)c\phi(x)}$, which maps $[0,\infty)$ onto $[0,\infty)$. 

The function $v$ is continuous, non-increasing and onto. We remark that Assumption 2 guarantees that $g$ (and thus $v$) is properly defined\footnote{Assumption 2 could in fact be relaxed by considering instead the inequality $(1-\rho)\phi_\infty<1/c$, which therefore enforces a constraint on the choice of both the $u$ function (through $\phi_\infty$) and the regularization parameter $\rho$. The proof of Theorem \ref{theorem1} (provided in Appendix) considers this more general case. Nevertheless, for simplicity of exposition, we will avoid this technical aspect in the core of the paper.}. Note that, importantly, $\phi_\infty$ does not have to be lower bounded by 1, as opposed to the non-regularized setting. With this in hand, we have the following theorem:
\begin{theorem} \label{theorem1}
Define $\mathcal{I}$ a compact set included in $ (0, 1]$. Let $\hat{\mathbf{C}}_N(\rho)$ be the unique solution to (\ref{implicit2}). Then, as $N,n \rightarrow \infty$, under Assumptions 1-2,
\[
\sup_{\rho \in \mathcal{I}}\left\lVert \hat{\mathbf{C}}_N(\rho)-\hat{\mathbf{S}}_N(\rho)\right\rVert \xrightarrow{\substack{\mathrm{a.s.}}} 0,
\]
where
\begin{align*}
\begin{array}{ll}
\hat{\mathbf{S}}_N(\rho) &\hspace{-2mm}\triangleq (1-\rho) v(\gamma)\frac{1}{n}\sum_{i=1}^{n}\mathbf{y}_i\mathbf{y}_i^\dagger+\rho \mathbf{I}_N, 
\end{array}
\end{align*}
with $\gamma$  the unique positive solution to the equation
\begin{align}
\label{eq:gamma_reg}
\begin{array}{rl}
\gamma&\hspace{-2mm}=\frac{1}{N}\tr \left[\mathbf{C}_N \left((1-\rho) \frac{v(\gamma)}{1+c (1-\rho) v(\gamma)\gamma}\mathbf{C}_N +\rho \mathbf{I}_N \right)^{-1}\right].
\end{array} 
\end{align}
Furthermore, the function $\rho \mapsto \gamma(\rho)$ is bounded, continuous on $(0,\infty]$ and greater than zero.
\end{theorem}
The proof of Theorem \ref{theorem1} (as well as that of the other technical results in this section) is provided in Appendix \ref{app:shrinkage}.
\begin{remark}
The uniform convergence in Theorem~\ref{theorem1} will be important for finding the optimal regularization parameter of a given estimator. As a matter of fact, the set $\mathcal{I}$, required to establish such uniform convergence, can be taken as $[0,1]$ in the over-sampled case (provided that $\liminf_N \{\lambda_1\left( \mathbf{C}_N \right) \} >0$).
\end{remark}

{\color{black}
\begin{remark}
Note that Theorem \ref{theorem1} (and Propositions \ref{cor 2}, \ref{cor 3} below) resemble similar results obtained in \cite{couillet2014large}. The key difference in the present work is that our results apply to a wide class of estimators (associated with $u$ functions satisfying the assumptions prescribed above), while \cite{couillet2014large} only focused on Tyler's estimator (associated with $u_\mathrm{Tyler}(x) = 1/x$). 
\end{remark}
}

Theorem \ref{theorem1} shows that, for every $u$ function, the estimator $\hat{\mathbf{C}}_N(\rho)$ asymptotically behaves (uniformly on $\rho \in  \mathcal{I}$) like the RSCM, with weights $\{(1-\rho) v(\gamma),\rho\}$ in lieu of the parameters $\{1-\beta,\beta\}$ in (\ref{LW}). Importantly, the relative weight given to the SCM depends on the underlying $u$ function, which entails that, for a fixed $\rho$, two different estimators may have different asymptotic behaviors.  However, while this is indeed the case, in the following it will be shown that, upon properly choosing the regularization parameter, all regularized M-estimators share the same, optimal asymptotic performance, at least with respect to a quadratic loss.

\subsection{Optimized regularization and asymptotic equivalence of different regularized M-estimators}
First, we will demonstrate that any trace-normalized regularized M-estimator is in fact asymptotically equivalent to the RSCM, up to a simple transformation of the regularization parameter. The result is as follows:

\begin{proposition}
\label{cor 2}
Let Assumptions 1-2 hold. For each $\rho \in (0,1]$, the parameter $\underline{\rho} \in (0,1]$ defined as 
\begin{align}
\label{eq:rho_rho_bar}
\underline{\rho} \triangleq \frac{\rho}{(1-\rho) v(\gamma)+\rho} 
\end{align} 
is such that 
\begin{align}
\label{eq:tn}
\frac{\hat{\mathbf{S}}_N(\rho)}{\frac{1}{N}\tr \hat{\mathbf{S}}_N(\rho)}=\frac{\mathbf{R}(\underline{\rho})}{\frac{1}{N}\tr \mathbf{R}(\underline{\rho})},
\end{align}
where we recall that $\mathbf{R}(\underline{\rho})=(1-\underline{\rho})\frac{1}{n}\sum_{i=1}^n \mathbf{y}_i \mathbf{y}_i^\dagger +\underline{\rho} \mathbf{I}_N$.

Reciprocally, for each $\underline{\rho} \in (0,1]$, there exists a solution $\rho \in (0,1]$ to the equation (\ref{eq:rho_rho_bar}) for which equality (\ref{eq:tn}) holds.
\end{proposition}

Proposition \ref{cor 2} implies that, in the absence of outliers, any (trace-normalized) estimator $\hat{\mathbf{S}}_N(\rho)$ is equal to a trace-normalized RSCM estimator with a regularization parameter $\underline{\rho}$ depending on $\rho$ and on the underlying $u$ function (through $v$). From Theorem \ref{theorem1}, it then follows that the estimator $\hat{\mathbf{C}}_N(\rho)$ asymptotically behaves like the RSCM estimator with parameter $\underline{\rho}$. 

Thanks to Proposition \ref{cor 2}, we may thus look for optimal asymptotic choices of $\rho$. Given an estimator $\hat{\mathbf{B}}_N$ of $\mathbf{C}_N$, define the quadratic loss of the associated trace-normalized estimator as:
\begin{align*}
\mathcal{L}\left(\frac{\hat{\mathbf{B}}_N}{\frac{1}{N}\tr\hat{\mathbf{B}}_N},\frac{\mathbf{C}_N}{\frac{1}{N}\tr \mathbf{C}_N}\right) \! \triangleq \! \frac{1}{N} \left\lVert \frac{\hat{\mathbf{B}}_N}{\frac{1}{N}\tr \hat{\mathbf{B}}_N}\!-\!\frac{\mathbf{C}_N}{\frac{1}{N}\tr \mathbf{C}_N}\right\rVert_F^2\!.
\end{align*}
We then have the following proposition:

\begin{proposition}
\label{cor 3}
\emph{(Optimal regularization)}

Let Assumptions 1 and 2 hold. Define
\begin{align*}
\mathcal{L}^\star &\triangleq c\frac{M_2 -1}{c+M_2-1}M_1^2 \\
\rho^\star &\triangleq \frac{c}{c+M_2-1},
\end{align*}
where $M_1 \triangleq \int t \nu(dt)$ and $M_2 \triangleq \int t^2 \nu(dt)$.
Then, 
\begin{align*}
\inf_{\rho \in \mathcal{I}} \mathcal{L}\left(\frac{\hat{\mathbf{C}}_N(\rho)}{\frac{1}{N}\tr\hat{\mathbf{C}}_N(\rho)},\frac{\mathbf{C}_N}{\frac{1}{N}\tr\mathbf{C}_N}\right) \xrightarrow{\substack{\mathrm{a.s.}}} \mathcal{L}^\star. 
\end{align*}
Furthermore, for $\hat{\rho}^\star$ a solution to $\frac{\hat{\rho}^\star}{(1-\hat{\rho}^\star)  v(\gamma)+\hat{\rho}^\star}=\rho^\star$, 
\[\mathcal{L}\left(\frac{\hat{\mathbf{C}}_N(\hat{\rho}^\star)}{\frac{1}{N}\tr\hat{\mathbf{C}}_N(\hat{\rho}^\star)},\frac{\mathbf{C}_N}{\frac{1}{N}\tr\mathbf{C}_N}\right) \xrightarrow{\substack{\mathrm{a.s.}}} \mathcal{L}^\star. \]
\emph{(Optimal regularization parameter estimate)}

The solution $\hat{\rho}_N \in \mathcal{I}$ to
\begin{align*}
\frac{\hat{\rho}_N}{\frac{1}{N}\tr \hat{\mathbf{C}}_N(\hat{\rho}_N)}=\frac{c_N}{\frac{1}{N} \tr \left[\left(\frac{1}{n}\sum_{i=1}^n \frac{\mathbf{y}_i\mathbf{y}_i^\dagger}{\frac{1}{N}\left\lVert\mathbf{y}_i\right\rVert^2} \right)^2\right]-1},
\end{align*}
satisfies
\begin{align*}
\hat{\rho}_N & \xrightarrow{\substack{\mathrm{a.s.}}} \hat{\rho}^\star \\ 
\mathcal{L}\left(\frac{\hat{\mathbf{C}}_N(\hat{\rho}_N)}{\frac{1}{N}\tr\hat{\mathbf{C}}_N(\hat{\rho}_N)},\frac{\mathbf{C}_N}{\frac{1}{N}\tr\mathbf{C}_N}\right) & \xrightarrow{\substack{\mathrm{a.s.}}}  \mathcal{L}^\star.
\end{align*}
\end{proposition}
Proposition \ref{cor 3} states that, irrespective of the choice of $u$, there exists some $\rho$ for which the quadratic loss of the corresponding regularized M-estimator is minimal, this minimum being the same as the minimum achieved by an optimally-regularized RSCM. The last result of Proposition~\ref{cor 3} provides a simple way to estimate this optimal parameter.

\begin{figure}[htb]
\begin{minipage}[b]{1\linewidth}
  \centering
  \centerline{\includegraphics[width=8cm]{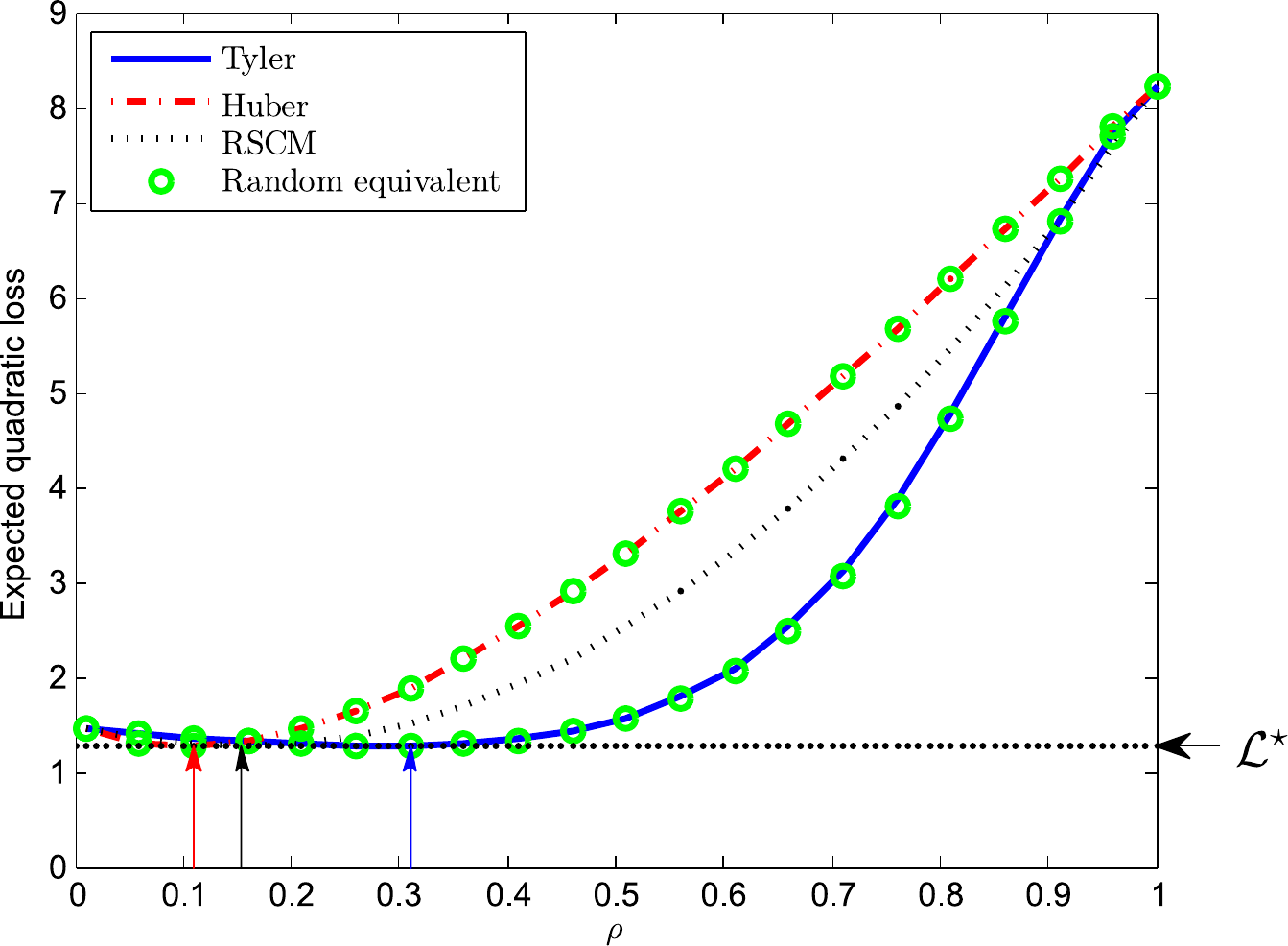}}
\end{minipage}
\caption{{\small Expected quadratic loss of different estimators as $\rho$ varies, for $N=150$, $n=100$ ($c_N=3/2$), and $[\mathbf{C}_N]_{ij}=.9^{|i-j|}$, averaged over 100 realizations. Arrows indicate the estimated optimal regularization parameters for the considered estimators, while $\mathcal{L}^\star$ indicates the asymptotic, minimal quadratic loss.}}
\label{fig:res10}
\end{figure}

In the following, we validate these theoretical findings through simulation. Let $[\mathbf{C}_N]_{ij}=.9^{|i-j|}$, and consider the $u$ functions specified in (\ref{u_tyler}) and (\ref{u_huber}), with $K=1/c_N$ and $t=0.1$.
For $\rho \in (0,1]$, Fig. \ref{fig:res10} depicts the expected quadratic loss $\mathcal{L}$ associated with the solution $\hat{\mathbf{C}}_N(\rho)$ of (\ref{implicit2}) and that associated with the RSCM (line curves), along with the expected quadratic loss associated with the random equivalent $\hat{\mathbf{S}}_N(\rho)$ of Tyler's and Huber's estimators (marker). 

For both $u$ functions and all $\rho \in (0,1]$, there is a close match between the quadratic loss of $\hat{\mathbf{C}}_N(\rho)$ and that of $\hat{\mathbf{S}}_N(\rho)$. This shows the accuracy of the (asymptotic) equivalence of $\hat{\mathbf{C}}_N(\rho)$ and $\hat{\mathbf{S}}_N(\rho)$ described in Theorem \ref{theorem1}. As suggested by our analysis, while the estimators associated with different $u$ functions have different performances for a given $\rho$, they have the same performance when $\rho$ is optimized, with a quadratic loss approaching $\mathcal{L}^\star$ for $N$ large. Furthermore, the optimal regularization parameter for a given $u$ function is accurately estimated, as shown by the arrows in Fig. \ref{fig:res10}.

\section{Large-dimensional robustness: Non-regularized case}
\label{sec:if}
In this section, we turn to the case where the data is contaminated by random outliers and study the robustness of M-estimators for distinct $u$ functions. Some initial insight has been previously provided in \cite{morales2015large} for the non-regularized case. Specifically, that study focused on the comparison of the weights given by the estimator to outlying and legitimate samples. Albeit insightful, the analysis in \cite{morales2015large} did not directly assess robustness, understood as the impact of outliers on the estimator's performance. {\color{black}Moreover, it is not clear how these weights translate into the Frobenius loss, or any other standard metric to measure the estimator's performance}.  Here we propose a different approach to analyze robustness, by introducing and evaluating a robustness metric which measures the bias induced by data contamination.  

We start by studying non-regularized estimators (thereby excluding the case $c_N \ge 1$), which are technically easier to handle. This will provide insight on the capabilities of different M-estimators to harness outlying samples. Then, in the following section, we will investigate how this study translates to the regularized case. The proofs of the technical results in this section are provided in Appendix \ref{app:if1}. 

\subsection{Asymptotic equivalent form under outlier data model}
We focus on a particular type of contamination model where outlying samples follow a distribution different from that of the legitimate samples. Similar to \cite{morales2015large}, the data matrix $\mathbf{Y}=[\mathbf{y}_1, \cdots, \mathbf{y}_{(1-\epsilon_n)n}, \mathbf{a}_1, \cdots, \mathbf{a}_{\epsilon_nn}] \in \mathbb{C}^N$ is constructed with the first $(1-\epsilon_n)n$ data samples ($\mathbf{y}_1, \cdots, \mathbf{y}_{(1-\epsilon_n)n}$) being the legitimate data, and following the same distribution as in Sections \ref{sec:2} and \ref{sec:shrinkage} (that is, they verify $\mathbf{y}_i=\mathbf{C}_N^{1/2}\mathbf{x}_i$). The remaining $\epsilon_nn$ ``contaminating'' data samples ($\mathbf{a}_1, \cdots, \mathbf{a}_{\epsilon_nn}$) are assumed to be random, independent of the $\mathbf{y}_i$, with $\mathbf{a}_i=\mathbf{D}_N^{1/2}\mathbf{x}_i'$, where $\mathbf{D}_N \in \mathbb{C}^{N\times N}$ is deterministic positive definite and $\mathbf{x}_1',\cdots,\mathbf{x}_{\epsilon_nn}'$ are independent random vectors with i.i.d.\@ zero mean, unit variance, and finite $(8+\eta)$-th order moment entries, for some $\eta > 0$. {\color{black} This outlier data model is effectively a cluster contamination model \cite{kaufman2009finding}.}

To characterize the asymptotic behavior of M-estimators for this data model, we require the following assumptions on the growth regime and on the underlying covariance matrices $\mathbf{C}_N$ and $\mathbf{D}_N$:
\begin{assumption}  
\hfill \begin{enumerate}[a.] 
\item $\epsilon_n \rightarrow \epsilon \in [0,1)$ and $c_N \rightarrow c \in (0, 1)$ as $N,n \rightarrow \infty$.

\item $ 0 < \liminf_N \{\lambda_1\left( \mathbf{C}_N \right) \} \leq \limsup_N \{\lambda_N\left( \mathbf{C}_N \right) \} < \infty$.
\item $\limsup_N \left\lVert \mathbf{D}_N \mathbf{C}_N^{-1}  \right\rVert < \infty$.
\end{enumerate}
\end{assumption}

Let us consider a function $u$ with now $1<\phi_{\infty}<\frac{1}{c}$. For such a $u$ function, the equation in $\mathbf{Z}$
\begin{align}
\mathbf{Z}&=\frac{1}{n}\sum_{i=1}^{(1-\epsilon_n)n}u\left(\frac{1}{N}\mathbf{y}_i^\dagger\mathbf{Z}^{-1}\mathbf{y}_i\right)\mathbf{y}_i\mathbf{y}_i^\dagger \nonumber \\
&+\frac{1}{n}\sum_{i=1}^{\epsilon_nn}u\left(\frac{1}{N}\mathbf{a}_i^\dagger\mathbf{Z}^{-1}\mathbf{a}_i\right)\mathbf{a}_i\mathbf{a}_i^\dagger \label{implicit_no_shrinkage}
\end{align}
has a unique solution \cite{maronna1976robust}, hereafter referred to as $\hat{\mathbf{C}}_N^{\epsilon_n}$. Moreover, define $v:[0,\infty) \rightarrow (0,u(0)]$ as $v(x)=u(g^{-1}(x))$ where $g^{-1}$ denotes the inverse function of $g(x)=\frac{x}{1-c\phi(x)}$, which maps $[0,\infty)$ onto $[0,\infty)$. 

In this setting, we have the following result:
\begin{theorem}  \label{theorem_no_shrinkage}
Let Assumption 3 hold and let $\hat{\mathbf{C}}_N^{\epsilon}$ be the unique solution to (\ref{implicit_no_shrinkage}). Then, as $N,n \rightarrow \infty$,
\[
\left\lVert \hat{\mathbf{C}}_N^{\epsilon_n}-\hat{\mathbf{S}}_N^{\epsilon_n}\right\rVert \xrightarrow{\substack{\mathrm{a.s.}}} 0
\]
where 
\begin{align*}
\hat{\mathbf{S}}_N^{\epsilon_n} \triangleq v(\gamma^{\epsilon_n})\frac{1}{n}\sum_{i=1}^{(1-\epsilon_n)n}\mathbf{y}_i\mathbf{y}_i^\dagger+v(\alpha^{\epsilon_n})\frac{1}{n}\sum_{i=1}^{\epsilon_nn}\mathbf{a}_i\mathbf{a}_i^\dagger, 
\end{align*}
with $\gamma^{\epsilon_n}$ and $\alpha^{\epsilon_n}$ the unique positive solutions to:
\begin{align*}
\gamma^{\epsilon_n}&=\frac{1}{N}\tr \mathbf{C}_N \mathbf{B}_N^{-1}\\
\alpha^{\epsilon_n}&=\frac{1}{N}\tr \mathbf{D}_N \mathbf{B}_N^{-1},
\end{align*}
with 
\begin{align*}
\mathbf{B}_N \triangleq \left(\frac{(1-\epsilon_n)v(\gamma^{\epsilon_n})}{1+cv(\gamma^{\epsilon_n})\gamma^{\epsilon_n}}\mathbf{C}_N +\frac{\epsilon_n v(\alpha^{\epsilon_n})}{1+cv(\alpha^{\epsilon_n})\alpha^{\epsilon_n}}\mathbf{D}_N\right).
\end{align*}

\end{theorem}
Theorem \ref{theorem_no_shrinkage} shows that $\hat{\mathbf{C}}_N^{\epsilon}$ behaves similar to a weighted SCM, with the legitimate samples weighted by $v(\gamma^{\epsilon_n})$, and the outlying samples by $v(\alpha^{\epsilon_n})$. 
This result generalizes \cite[Corollary 3]{morales2015large} to allow for $\epsilon \in [0,1)$. {\color{black} Indeed, the applicability of the results of \cite{morales2015large} was limited by a constraint on $\epsilon, \phi_\infty$ and $c$, which had to verify the inequality $(1-\epsilon)^{-1}<\phi_{\infty}<\frac{1}{c}$ (along with $c <1-\epsilon$). This constraint  can be easily violated in practice, particularly for $c \lesssim 1$: 
for example, if $c = 0.9$, then having a proportion of outliers as little as $15\%$ 
would imply that $\epsilon  > 1-c$, in which case there does not exist a $\phi_{\infty}$ verifying the inequality above.

A scenario which will be of particular interest in the following concerns the case where there is a vanishingly small proportion of outliers. This occurs when $\epsilon_n = O(1/n^\mu)$ for some $0<\mu \leq 1$, in which case $\epsilon_n \rightarrow \epsilon= 0$. For this scenario, the weights given to the legitimate and outlying data are
\begin{align}
\label{eq:limit_gamma1}
\gamma^0 \triangleq & \lim_{n \rightarrow \infty} \gamma^{\epsilon_n}  = \frac{\phi^{-1}(1)}{1-c}  \\
\label{eq:limit_alpha1}
\alpha^{0}  \triangleq & \lim_{n \rightarrow \infty} \alpha^{\epsilon_n}  = \gamma^{0} \frac{1}{N} \tr (\mathbf{C}_N^{-1}\mathbf{D}_N),
\end{align}
respectively.

In the following, we exploit the form of $\hat{\mathbf{S}}_N^{\epsilon_n}$ to characterize the effect of random outliers on the estimator $\hat{\mathbf{C}}_N^{\epsilon_n}$. 

\subsection{Robustness analysis}
Let $\hat{\mathbf{C}}_N^0$ be the solution to (\ref{eq:implicit Z}), and $\hat{\mathbf{C}}_N^{\epsilon_n}$ the solution to (\ref{implicit_no_shrinkage}). 
We propose the following metric, termed \emph{measure of influence}, to assess the robustness of a given estimator to an $\epsilon$-contamination of the data:\\

\begin{definition}
For $\epsilon_n \rightarrow \epsilon \in [0,1)$, the measure of influence $\mathrm{MI}(\epsilon_n)$ is given by
\begin{align*}
\mathrm{MI}(\epsilon_n) \triangleq \left\lVert \mathbb{E}\left[\frac{\hat{\mathbf{C}}_N^{0}}{\frac{1}{N}\tr \hat{\mathbf{C}}_N^{0}}-\frac{\hat{\mathbf{C}}_N^{\epsilon_n}}{\frac{1}{N}\tr \hat{\mathbf{C}}_N^{\epsilon_n}}\right] \right\rVert.
\end{align*}
\end{definition}

For simplicity, we assume hereafter that $\frac{1}{N}\tr \mathbf{C}_N = \frac{1}{N}\tr \mathbf{D}_N = 1 $ for all $N$. From Theorems \ref{th:Romain_no_shrinkage} and \ref{theorem_no_shrinkage}, we have the following:
\begin{corollary}
\label{cor:mi_nonreg}
As $N,n \rightarrow \infty$,
\begin{align*}
\mathrm{MI}(\epsilon_n)-\overline{\mathrm{MI}}(\epsilon_n)\rightarrow 0,
\end{align*}
where
\begin{align}
\overline{\mathrm{MI}}(\epsilon_n) = \frac{\epsilon_n v(\alpha^{\epsilon_n})}{(1-\epsilon_n)v(\gamma^{\epsilon_n})+\epsilon_n v(\alpha^{\epsilon_n})} \left\lVert \mathbf{C}_N-\mathbf{D}_N  \right\rVert.
\label{eq:F}
\end{align}
\end{corollary}

Note that $\lim_{\epsilon_n \rightarrow 0} \overline{\mathrm{MI}}(\epsilon_n) = \overline{\mathrm{MI}}(0) = 0$, as expected. The result (\ref{eq:F}) shows that $\overline{\mathrm{MI}}$ is globally influenced by $\left\lVert \mathbf{C}_N-\mathbf{D}_N  \right\rVert$, which is also an intuitive result, since it suggests that the more ``different'' $\mathbf{D}_N$ is from $\mathbf{C}_N$, the higher the influence of the outliers on the estimator. To get clearer insight on the effect of a small proportion of outliers, assuming $\epsilon_n = O(1/n^\mu)$ for some $0<\mu \leq 1$, we compute
\begin{align}
\label{def:ia}
\overline{\mathrm{IMI}} \triangleq \lim_{n \rightarrow \infty} \frac{1}{\epsilon_n} \overline{\mathrm{MI}}(\epsilon_n),
\end{align} 
which we will refer to as the infinitesimal measure of influence (IMI)\footnote{{\color{black}As defined, the IMI is reminiscent of a standard robust statistical metric known as the influence function \cite{hampel2011robust}. A key distinction with our proposed metric is that the influence function measures the effect of an infinitesimal mass contamination $\{\mathbf{x}\}$ on the estimator robustness (here, $\mathbf{x}$ is an arbitrary point). In contrast, the IMI focuses on the impact of an infinitesimal \emph{random} contamination.}}.

From (\ref{eq:F}) and (\ref{def:ia}),
\begin{align}
\label{eq:if_non_reg}
\overline{\mathrm{IMI}} = \frac{v(\alpha^0)}{v(\gamma^0)} \left\lVert \mathbf{C}_N-\mathbf{D}_N  \right\rVert,
\end{align}
with $\gamma^0$, $\alpha^0$ given in (\ref{eq:limit_gamma1}) and (\ref{eq:limit_alpha1}), respectively.

For particular $u$ functions, these general results reduce to even simpler forms: for example, for $u$ functions such that $\phi^{-1}(1) = 1$ \big(such as $u_\mathrm{M-Tyler}=\frac{1+t}{t+x}$ or $u_\mathrm{M-Huber} = \min\{1,\frac{1+t}{t+x}\}$\big), which entails $v(\gamma^0) = 1$, (\ref{eq:if_non_reg}) further yields 
\begin{align*}
\overline{\mathrm{IMI}} 
 = v \left( \frac{\frac{1}{N} \tr \mathbf{C}_N^{-1}\mathbf{D}_N}{1-c}\right) \left\lVert \mathbf{C}_N-\mathbf{D}_N  \right\rVert.
 \end{align*}
Further, considering $t$ small, the IMI associated with $u_\mathrm{M-Tyler}$ and $u_\mathrm{M-Huber}$ can be approximated as
 \begin{align}
\label{eq:if_T}
\overline{\mathrm{IMI}}_\mathrm{M-Tyler} & \simeq \frac{1+t}{t+\frac{1}{N}\tr \mathbf{C}_N^{-1} \mathbf{D}_N} \left\lVert \mathbf{C}_N-\mathbf{D}_N \right\rVert
\end{align}
and
 \begin{align}
 \label{eq:if_H}
\overline{\mathrm{IMI}}_\mathrm{M-Huber} & \simeq \left\{ \begin{array}{ll}
  \left\lVert \mathbf{C}_N-\mathbf{D}_N \right\rVert & \text{if } \frac{1}{N}\tr \mathbf{C}_N^{-1} \mathbf{D}_N \leq 1 \\
  \overline{\mathrm{IMI}}_\mathrm{M-Tyler} & \text{if } \frac{1}{N}\tr \mathbf{C}_N^{-1} \mathbf{D}_N > 1
  \end{array}  \right..
\end{align}
Hence, when $\frac{1}{N}\tr \mathbf{C}_N^{-1} \mathbf{D}_N \leq 1$, $\overline{\mathrm{IMI}}_\mathrm{M-Huber}  \leq \overline{\mathrm{IMI}}_\mathrm{M-Tyler}$, which shows that the influence of an infinitesimal fraction of outliers is higher for Tyler's estimator than for Huber's. In contrast, when $\frac{1}{N}\tr \mathbf{C}_N^{-1} \mathbf{D}_N > 1$, both Huber's and Tyler's estimators exhibit the same IMI.

For comparison, the measure of influence of the SCM can be written as
\begin{align*}
\overline{\mathrm{MI}}_\mathrm{SCM}(\epsilon_n) = \epsilon_n \left\lVert \mathbf{C}_N-\mathbf{D}_N  \right\rVert,
\end{align*}
which is linear in $\epsilon_n$. 
It follows immediately that
\begin{align}
\overline{\mathrm{IMI}}_\mathrm{SCM} = \left\lVert \mathbf{C}_N-\mathbf{D}_N  \right\rVert.
\label{eq:ia_scm}
\end{align}
The fact that $\overline{\mathrm{IMI}}_\mathrm{SCM}$ is bounded may seem surprising, since it is known that a single arbitrary outlier can arbitrarily bias the SCM \cite{hampel2011robust}, however we recall that the current model focuses on a particular random outlier scenario.
From (\ref{eq:F}), the SCM is more affected than given M-estimators by the introduction of outliers if and only if 
\begin{align*}
\overline{\mathrm{MI}}(\epsilon_n) \leq \overline{\mathrm{MI}}_\mathrm{SCM}(\epsilon_n)  \Leftrightarrow v(\alpha^{\epsilon_n}) \leq v(\gamma^{\epsilon_n}).
\end{align*}
This further legitimizes the study in \cite{morales2015large}, which focused on these weights to assess the robustness of a given M-estimator. However, in the regularized case it will be shown that the relationship between the relative weights and robustness is more complex (see Subsection \ref{sec:rob_reg}). 

Fig.\@ \ref{fig:F_full_rank} depicts the measure of influence $\overline{\mathrm{MI}}(\epsilon_n)$ for different $u$ functions, as the proportion $\epsilon_n$ of outlying samples increases. For every $u$ function ($u_\mathrm{M-Tyler}$ or $u_\mathrm{M-Huber}$), we take $t = 0.1$. In addition, we show the measure of influence of the SCM, as well as the linear approximation $\epsilon_n \mapsto \epsilon_n \overline{\mathrm{IMI}}$ (computed using (\ref{eq:if_T}), (\ref{eq:if_H}) and (\ref{eq:ia_scm})) of the measure of influence in the neighborhood of $\epsilon = 0$. We first set $[\mathbf{C}_N]_{ij}=.9^{|i-j|}$ and $[\mathbf{D}_N]_{ij}=.2^{|i-j|}$ (such that $\frac{1}{N}\tr \mathbf{C}_N^{-1} \mathbf{D}_N > 1$), and then swap the roles of $\mathbf{C}_N$ and $\mathbf{D}_N$ (such that $\frac{1}{N}\tr \mathbf{C}_N^{-1} \mathbf{D}_N < 1$). 

\begin{figure*}[h]
\centering
\begin{subfigure}{.45\textwidth}
  \centering
  \includegraphics[width=.95\linewidth]{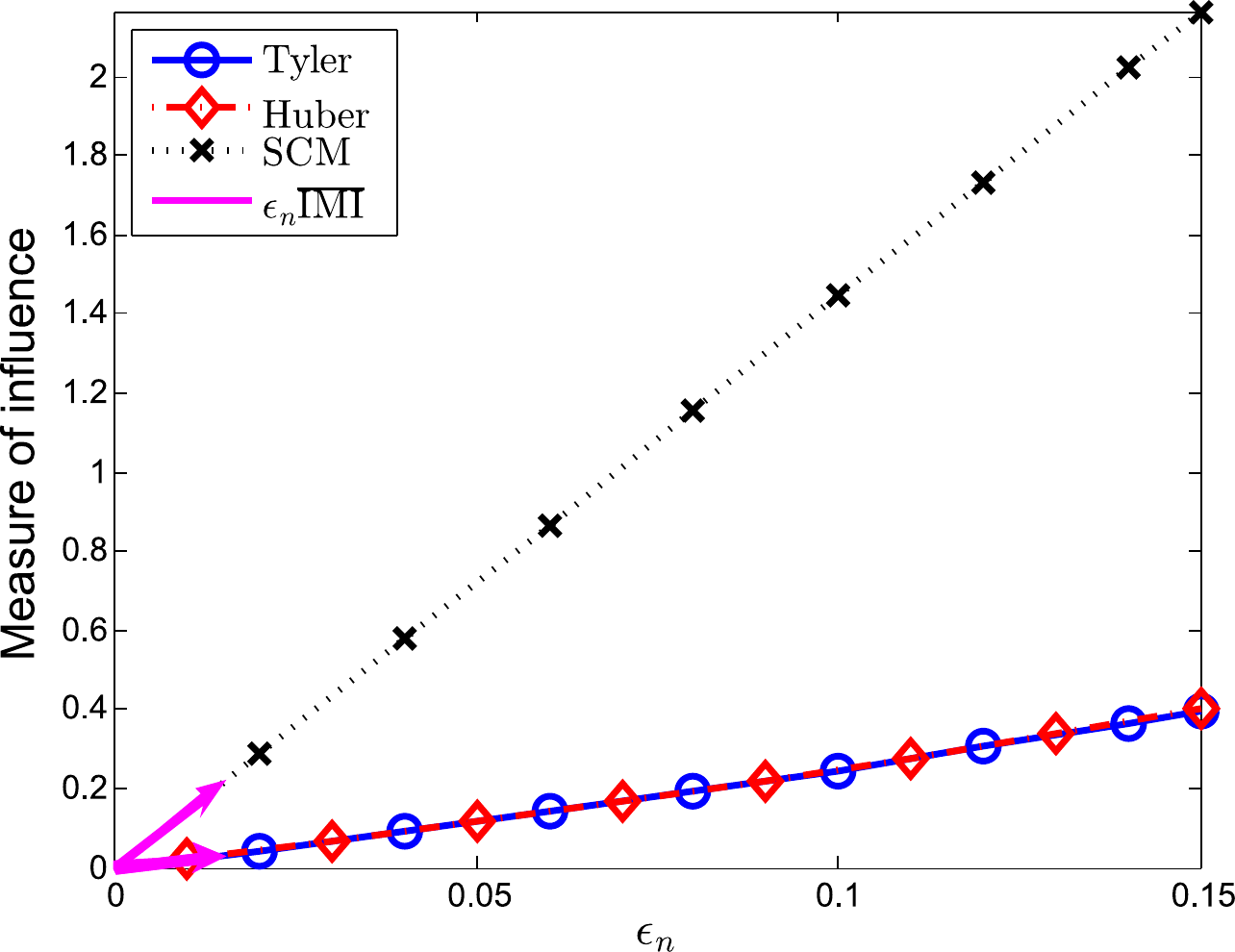}
  \caption{$[\mathbf{C}_N]_{ij}=.9^{|i-j|}$ and $[\mathbf{D}_N]_{ij}=.2^{|i-j|}$.}
  \label{fig:sub1}
\end{subfigure}%
\begin{subfigure}{.45\textwidth}
  \centering
  \includegraphics[width=.95\linewidth]{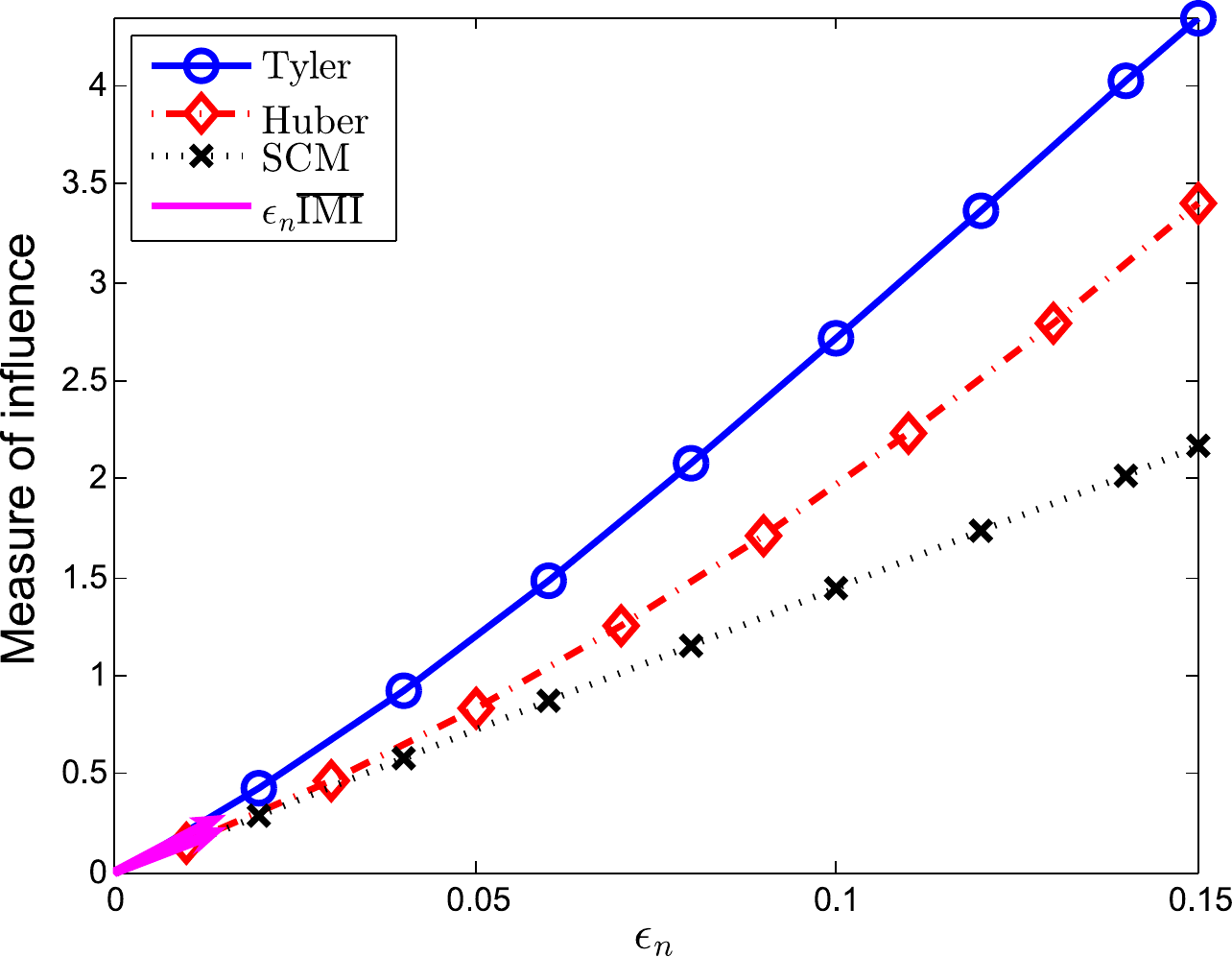}
  \caption{$[\mathbf{C}_N]_{ij}=.2^{|i-j|}$ and $[\mathbf{D}_N]_{ij}=.9^{|i-j|}$.}
  \label{fig:sub2}
\end{subfigure}
\caption{{\small Measure of influence for $\epsilon_n \in [0, 0.15]$, in the non-regularized case for $N=50$, $n=200$ ($c_N = 1/4$).}}
\label{fig:F_full_rank}
\end{figure*}

In the case where $\frac{1}{N}\tr \mathbf{C}_N^{-1} \mathbf{D}_N > 1$, Fig.\@ \ref{fig:F_full_rank} confirms that the measure of influence of both Tyler's and Huber's estimators is lower than that of the SCM, as corroborated by the fact that  $\overline{\mathrm{IMI}}_\mathrm{M-Tyler},\overline{\mathrm{IMI}}_\mathrm{M-Huber} < \left\lVert \mathbf{C}_N-\mathbf{D}_N  \right\rVert = \overline{\mathrm{IMI}}_\mathrm{SCM}$ (see (\ref{eq:if_T}), (\ref{eq:if_H})). This shows that the considered M-estimators are more robust to the introduction of outliers than the SCM. Furthermore both Tyler's and Huber's estimators exhibit the same robustness for small $\epsilon_n$. However, in the opposite case where $\frac{1}{N}\tr \mathbf{C}_N^{-1} \mathbf{D}_N <1$, Tyler's estimator is much less robust than both Huber's estimator and the SCM, which are both equally robust (for small $\epsilon_n$). Since both $\mathbf{C}_N$ and $\mathbf{D}_N$ are unknown in practice, it suggests that choosing Huber's estimator is preferable over Tyler's. 

\section{Large dimensional robustness: Regularized case}
\label{sec:if_reg}
We now turn to the regularized case, which, in particular, allows $c \ge 1$. The proofs of the technical results in this section are provided in Appendix \ref{app:if2}.

\vspace{-0.1cm}
\subsection{Asymptotic equivalent form under outlier data model}
\label{sec:out}
To facilitate the robustness study of regularized M-estimators, we start by analyzing the large-dimensional behavior of these estimators in the presence of outliers.

For $\rho \in \mathcal{R} = (\rho_0,1]$, where $\rho_0=\max{\{0,1-\frac{1}{c\phi_\infty}\}}$, we define the regularized estimator $\hat{\mathbf{C}}^\epsilon_N(\rho)$ associated with the function $u$ as the unique solution to the equation in $\mathbf{Z}$:
\begin{align}
\mathbf{Z}&=(1-\rho)\frac{1}{n}\sum_{i=1}^{(1-\epsilon_n)n}u\left(\frac{1}{N}\mathbf{y}_i^\dagger\mathbf{Z}^{-1}\mathbf{y}_i\right)\mathbf{y}_i\mathbf{y}_i^\dagger \nonumber\\
&+(1-\rho)\frac{1}{n}\sum_{i=1}^{\epsilon_nn}u\left(\frac{1}{N}\mathbf{a}_i^\dagger\mathbf{Z}^{-1}\mathbf{a}_i\right)\mathbf{a}_i\mathbf{a}_i^\dagger+\rho\mathbf{I}_N. \label{implicit_rho_version}
\end{align}

\begin{remark}
If we assume that $\phi_\infty < \frac{1}{c}$, then the range of admissible $\rho$ is $\mathcal{R} =(0,1]$. Furthermore, if $c < 1$, we can in fact take $\mathcal{R}=[0,1]$. In the following, we assume that $\phi_\infty < \frac{1}{c}$. Note that, similar to the outlier-free scenario, the introduction of a regularization parameter allows us to relax the assumption of $\phi_\infty>1$.
\end{remark}

\begin{theorem}  \label{theorem2}
Assume the same contaminated data model as in Theorem \ref{theorem_no_shrinkage}. Let Assumptions 1-2 hold and let $\hat{\mathbf{C}}_N^{\epsilon_n}(\rho)$ be the unique solution to (\ref{implicit_rho_version}). Then, as $N,n \rightarrow \infty$, for all $\rho \in \mathcal{R}$,
\[
\left\lVert \hat{\mathbf{C}}_N^{\epsilon_n}(\rho)-\hat{\mathbf{S}}_N^{\epsilon_n}(\rho)\right\rVert \xrightarrow{\substack{\mathrm{a.s.}}} 0
\]
where 
\begin{align*}
\hat{\mathbf{S}}_N^{\epsilon_n}(\rho) &\triangleq (1-\rho)v(\gamma^{\epsilon_n})\frac{1}{n}\sum_{i=1}^{(1-\epsilon_n)n}\mathbf{y}_i\mathbf{y}_i^\dagger \\
&+(1-\rho)v(\alpha^{\epsilon_n})\frac{1}{n}\sum_{i=1}^{\epsilon_nn}\mathbf{a}_i\mathbf{a}_i^\dagger+\rho\mathbf{I}_N, 
\end{align*}
with $\gamma^{\epsilon_n}$ and $\alpha^{\epsilon_n}$ the unique positive solutions to:
\begin{align}
\label{eq:rd_shrinkage1}
\gamma^{\epsilon_n}&=\frac{1}{N}\tr \mathbf{C}_N \mathbf{B}_N^{-1} \\
\alpha^{\epsilon_n}&=\frac{1}{N}\tr \mathbf{D}_N \mathbf{B}_N^{-1}, \nonumber
\end{align}
with 
\begin{align*}
\mathbf{B}_N &\triangleq (1-\rho)\frac{(1-\epsilon_n)v(\gamma^{\epsilon_n})}{1+(1-\rho)cv(\gamma^{\epsilon_n})\gamma^{\epsilon_n}}\mathbf{C}_N \\
&+(1-\rho)\frac{\epsilon_n v(\alpha^{\epsilon_n})}{1+(1-\rho)cv(\alpha^{\epsilon_n})\alpha^{\epsilon_n}}\mathbf{D}_N+\rho\mathbf{I}_N.
\end{align*}
\end{theorem}

\begin{remark}
In the case $\epsilon_n=0$ (no outliers), Theorem \ref{theorem2} reduces to Theorem \ref{theorem1}, while in the case $\rho = 0$ (if $c < 1$), it reduces to Theorem \ref{theorem_no_shrinkage}.
\end{remark}

\subsection{Robustness analysis}
\label{sec:rob_reg}
Similar to the non-regularized case, we next make use of $\hat{\mathbf{S}}_N^{\epsilon_n}(\rho)$ to study the robustness of $\hat{\mathbf{C}}_N^{\epsilon_n}(\rho)$. Importantly, introducing a regularization parameter entails an additional variable to consider when studying the robustness of M-estimators.

We denote by $\hat{\mathbf{C}}_N^0(\rho)$ the solution to (\ref{implicit_rho_version}) for a given $\rho \in  \mathcal{R}$ when there are no outliers. Similar to the non-regularized case, we define
\begin{align*}
\mathrm{MI}(\rho,\epsilon_n) \triangleq \left\lVert \mathbb{E}\left[ \frac{\hat{\mathbf{C}}_N^0(\rho)}{\frac{1}{N}\tr\hat{\mathbf{C}}_N^0(\rho)}-\frac{\hat{\mathbf{C}}_N^{\epsilon_n}(\rho)}{\frac{1}{N}\tr\hat{\mathbf{C}}_N^{\epsilon_n}(\rho)}\right] \right\rVert.
\end{align*}

By Theorem \ref{theorem2}, we have the following corollary:
\begin{corollary}
\label{cor:mi_shrinkage}
Let the same assumptions as in Theorem \ref{theorem2} hold. Then,
\begin{align*}
\forall \rho \in \mathcal{R}, \quad \mathrm{MI}(\rho,\epsilon_n) - \overline{\mathrm{MI}}(\rho,\epsilon_n) \rightarrow 0,
\end{align*}
with
\begin{align*}
\overline{\mathrm{MI}}(\rho,\epsilon_n)  = \left\lVert\frac{U(\epsilon_n,\rho)}{V(\epsilon_n,\rho)}\right\rVert,
\end{align*}
with $U(\epsilon_n,\rho),V(\epsilon_n,\rho)$ defined as
\begin{align*}
U(\epsilon_n,\rho) &= \rho (1-\rho) ((1-\epsilon_n)v(\gamma^{\epsilon_n})-v(\gamma^0))(\mathbf{C}_N-\mathbf{I}_N)\\
&+\rho (1-\rho) \epsilon_n v(\alpha^{\epsilon_n}) (\mathbf{D}_N-\mathbf{I}_N) \\
&+ (1-\rho)^2 \epsilon_n v(\gamma^0) v(\alpha^{\epsilon_n}) (\mathbf{D}_N-\mathbf{C}_N)\\
V(\epsilon_n,\rho) &= ((1-\rho)(1-\epsilon_n)v(\gamma^{\epsilon_n})\\
&+(1-\rho)\epsilon_n v(\alpha^{\epsilon_n}) + \rho)((1-\rho)v(\gamma^0)+\rho).
\end{align*}
\end{corollary}

Unlike in the non-regularized case, the form of $\overline{\mathrm{MI}}(\rho,\epsilon_n)$ renders the analysis difficult in general. For the specific case $\rho = 1$ however, $\overline{\mathrm{MI}}(1,\epsilon_n) = 0$ for all $\epsilon_n$. This is intuitive, and reflects the fact that the more we regularize an estimator, the more robust it becomes (eventually, it boils down to taking $\hat{\mathbf{C}}_N = \mathbf{I}_N$). This extreme regularization, however, leads to a significant bias, and is therefore not desirable.

In the following, we focus on the infinitesimal measure of influence associated with $\mathrm{MI}(\rho,\epsilon_n)$, which is defined in a similar way as for the non-regularized case: assume $\epsilon_n = O(1/n^\mu)$ for some $0<\mu \leq 1$ ($\epsilon_n \rightarrow \epsilon = 0$). Then, for $v$ smooth enough\footnote{Precise details are provided in Appendix \ref{app:if}.}, and $\rho \in \mathcal{R}$,
\begin{align*}
\overline{\mathrm{IMI}}(\rho) \triangleq \lim_{n \rightarrow \infty} \frac{1}{\epsilon_n} \overline{\mathrm{MI}}(\rho,\epsilon_n).
\end{align*}

\noindent Corollary \ref{cor:mi_shrinkage} allows us to compute $\overline{\mathrm{IMI}}(\rho)$ explicitly in the particular case where $\frac{1}{N}\tr \mathbf{C}_N = \frac{1}{N}\tr \mathbf{D}_N = 1 $ for all $N$. This is given as follows:
\begin{corollary} Let the same assumptions as in Theorem \ref{theorem2} hold. 
If $\gamma \mapsto v(\gamma)$ is differentiable in the neighborhood of $\gamma^0 = \lim_{n \rightarrow \infty} \gamma^{\epsilon_n}$,
\label{cor:ia_shrinkage} 
\begin{align}
\overline{\mathrm{IMI}}(\rho) = \frac{1}{((1-\rho)v(\gamma^0)+\rho)^2}\left\lVert \mathbf{G}(\rho) \right\rVert, 
\label{eq:if_shrinkage}
\end{align}
where
\begin{align*}
\mathbf{G}(\rho) &= \rho(1-\rho) [v(\alpha^0) (\mathbf{D}_N-\mathbf{I}_N) - v(\gamma^0) (\mathbf{C}_N-\mathbf{I}_N)] \\
& + (1-\rho)^2v(\gamma^0)v(\alpha^0)(\mathbf{D}_N-\mathbf{C}_N) \\
& +  \rho(1-\rho) \left. \frac{dv}{d\epsilon_n} \right|_{\epsilon_n = 0}(\mathbf{C}_N-\mathbf{I}_N),
\end{align*}
and where 
\begin{align*}
\alpha^0 & \triangleq \lim_{n \rightarrow \infty} \alpha^{\epsilon_n} \\
&= \frac{1}{N} \tr \mathbf{D}_N \left(\frac{(1-\rho)v(\gamma^0)}{1+(1-\rho)c\gamma^0v(\gamma^0)}\mathbf{C}_N+\rho \mathbf{I}_N \right)^{-1}.
\end{align*}
\end{corollary}

Details on how to evaluate $\left.\frac{dv}{d\epsilon} \right|_{\epsilon = 0}$ are given in Appendix~\ref{app:if}. 

While the intricate expression $\mathbf{G}(\rho)$ does not yield simple analytical insight for an arbitrary regularization parameter $\rho \in \mathcal{R}$, it can still be leveraged to numerically assess the robustness of regularized estimators, as we show below. As a given $\rho$ plays an a-priori different role for distinct M-estimators, a direct comparison of $\overline{\mathrm{MI}}(\rho, \epsilon_n)$ or $\overline{\mathrm{IMI}}(\rho)$ (for fixed $\rho$) for different estimators  is not meaningful. However, using Proposition \ref{cor 3}, we can choose $\rho = \rho^\star$ such that, in the absence of outliers, a given estimator's quadratic loss is minimal. This allows us to meaningfully compare how robust these estimators are to the introduction of a small proportion of outliers. 

For our subsequent numerical studies, we will focus on the scenario where $\frac{1}{N} \tr \mathbf{C}_N^{-1}\mathbf{D}_N > 1$. In the alternative case ($\frac{1}{N} \tr \mathbf{C}_N^{-1}\mathbf{D}_N < 1$), the differences between Huber, Tyler, and the RSCM are marginal (at least for small $\epsilon$). In Fig.\@~\ref{fig:if_reg0}, we compute the measure of influence $\overline{\mathrm{MI}}(\rho^\star,\epsilon_n)$ of the RSCM and of the estimators associated with $u_\mathrm{M-Tyler}$ and $u_\mathrm{M-Huber}$ (with $K=1/c_N$), as $\epsilon_n$ increases. We also plot the linear approximation $\epsilon_n \mapsto \epsilon_n \overline{\mathrm{IMI}}(\rho^\star)$ (computed using (\ref{eq:if_T}), (\ref{eq:if_H}) and (\ref{eq:ia_scm})) of $\overline{\mathrm{MI}}$ in the neighborhood of $\epsilon_n = 0$. We observe that the MI of Tyler's estimator is lower than that of Huber's estimator. This differs from the non-regularized case, where Tyler's and Huber's IMI were shown to be the same. This suggests that ``less-correlated'' outlying samples have a greater negative impact on regularized Huber's estimator, as compared with Tyler's estimator. It also appears that $\epsilon_n \overline{\mathrm{IMI}}(\rho^\star)$ is a fairly good approximation of $\overline{\mathrm{MI}}(\rho^\star,\epsilon_n)$ for small $\epsilon_n$, which shows the interest of Corollary \ref{cor:ia_shrinkage}. 
\begin{figure}[!]
\centering
  \includegraphics[width=8cm]{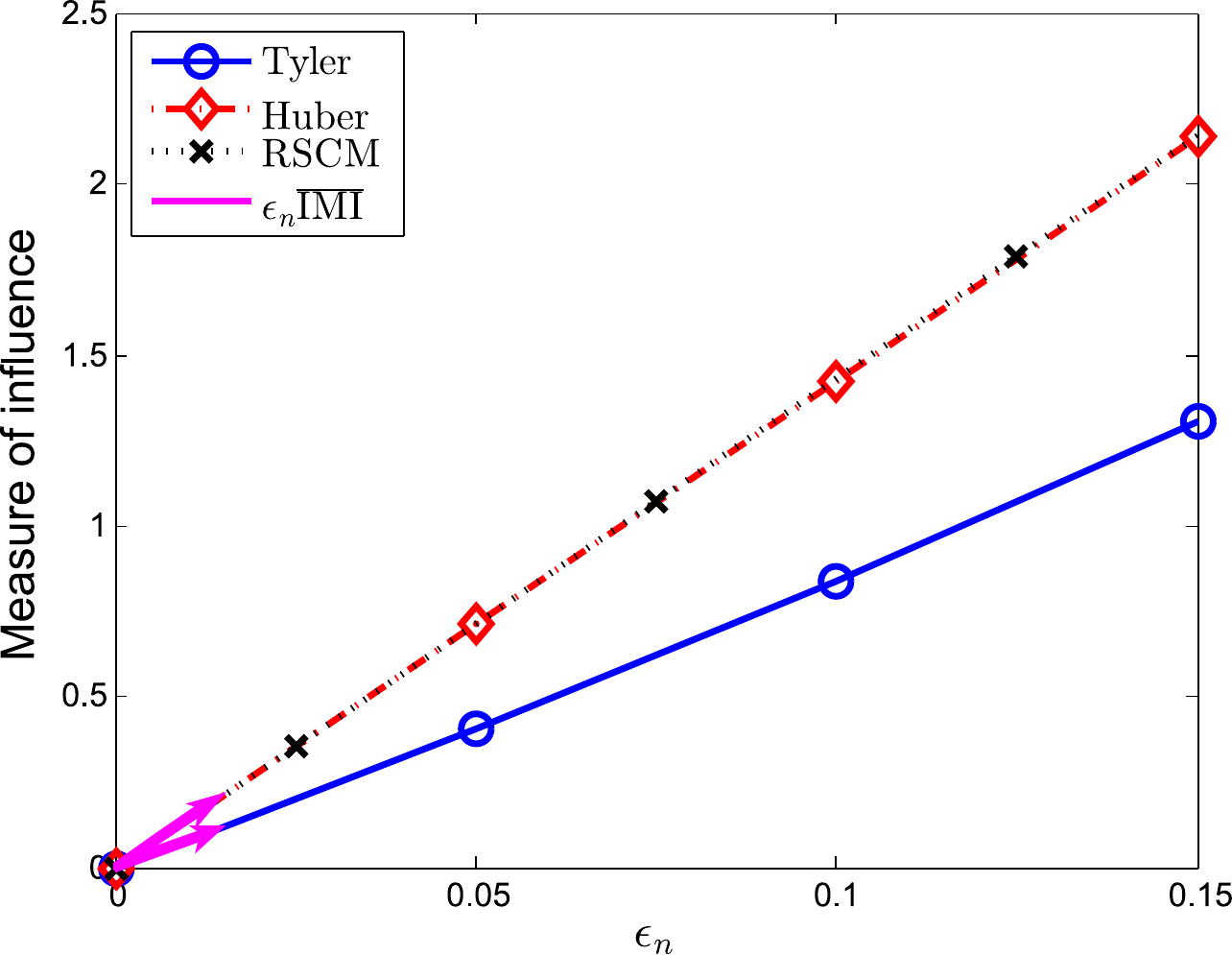}
\caption{{\small Measure of influence for $\epsilon_n \in [0, 0.15]$, in the regularized case ($N=150$, $n=100$, such that $c_N = 3/2$). The IMI is computed at the optimal regularization parameter (assuming clean data) that minimizes the quadratic loss of the estimator. $[\mathbf{C}_N]_{ij}=.9^{|i-j|}$ and $[\mathbf{D}_N]_{ij}=.2^{|i-j|}$.}}
\label{fig:if_reg0}
\end{figure}

So far, we have only considered two possible values of $c_N$: $c_N = 1/4$ (non-regularized case, Fig.\@ \ref{fig:F_full_rank}) and $c_N = 3/2$ (regularized case, Fig.\@ \ref{fig:if_reg0}).
To connect our results in the regularized and non-regularized scenarios, we now  evaluate $\overline{\mathrm{IMI}}(\rho^\star)$ for various $c_N$. Such experiment shall shed light on the effect of the aspect ratio $c_N$ on the robustness of different estimators. We consider again $u_\mathrm{M-Tyler}$ and $u_\mathrm{M-Huber}$, but now with $K=\min \{1,\frac{1}{c_N}\}$, such that Assumption 2 is verified; note that for $c_N \leq 1$, we retrieve the setting of Fig.\@ \ref{fig:F_full_rank}. Results are reported in Fig.\@ \ref{fig:if_reg2}. It appears that the IMI of a regularized estimator varies with $c_N$ in a non-trivial manner. Indeed, different $c_N$ call for different amounts of regularization (through $\rho^\star$), which in turn lead to substantial differences in robustness. Nonetheless, when $c_N \rightarrow 0$, the IMI of a given estimator tends to its non-regularized counterpart (indicated by arrows). This is a natural result, since in such case the need for regularization vanishes. We also notice that Tyler's estimator shows better robustness than all other estimators for nearly all values of $c_N$. 

\begin{figure}[h]
\centering
  \includegraphics[width=8cm]{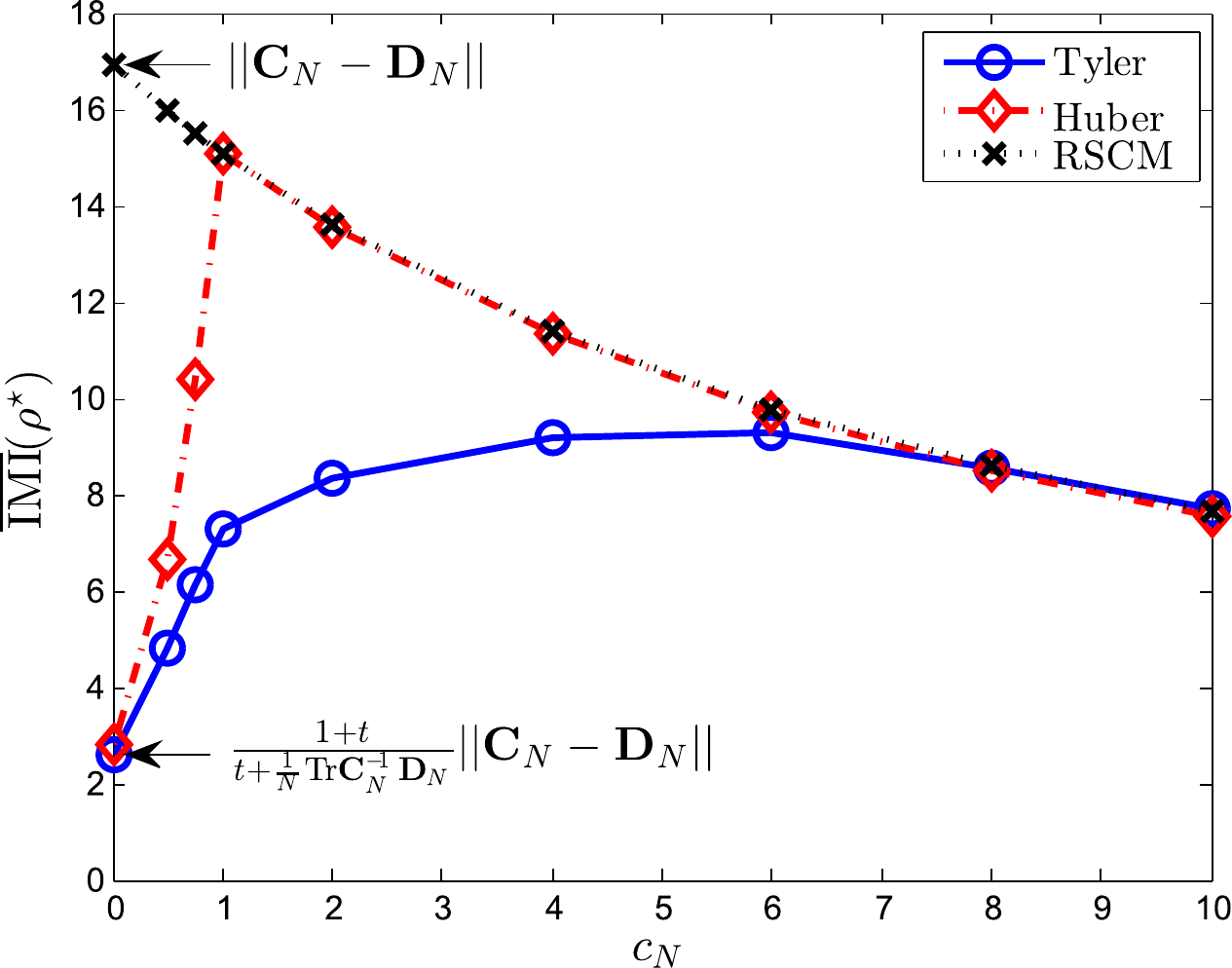}
\caption{{\small Infinitesimal measure of influence vs.\@ the aspect ratio $c_N$. The IMI is computed at the optimal regularization parameter (assuming clean data) that minimizes the quadratic loss of the estimator. Arrows indicate the IMI in the non-regularized case ($\rho =0$). $[\mathbf{C}_N]_{ij}=.9^{|i-j|}$ and $[\mathbf{D}_N]_{ij}=.2^{|i-j|}$.}}
\label{fig:if_reg2}
\end{figure}

\section{Discussion and concluding remarks} 
In summary, we have shown that, in the absence of outliers, regularized M-estimators are asymptotically equivalent to RSCM estimators and that, when optimally regularized, they all attain the same performance as the optimal RSCM, at least with respect to the quadratic loss. We proposed an intuitive metric to assess the robustness of different estimators when random outliers are introduced. In particular, it was shown in the non-regularized case that Huber's estimator is generally preferable over Tyler's, while in the regularized case,  {\color{black} if estimators are optimally regularized, Tyler is preferable over Huber when legitimate samples are ``more correlated'' than outlying samples.}

The comparatively different behaviour in regularized and non-regularized settings evidences the substantial (and non-trivial) effect of regularization on the robustness of M-estimators. This point is further emphasized in Fig.\@ \ref{fig:if_reg1}, where 
we plot $\overline{\mathrm{IMI}}(\rho)$ for $\rho \in (0,1]$.
\begin{figure}[htb]
  \begin{minipage}[b]{1\linewidth}
  \centering
  \centerline{\includegraphics[width=7.5cm]{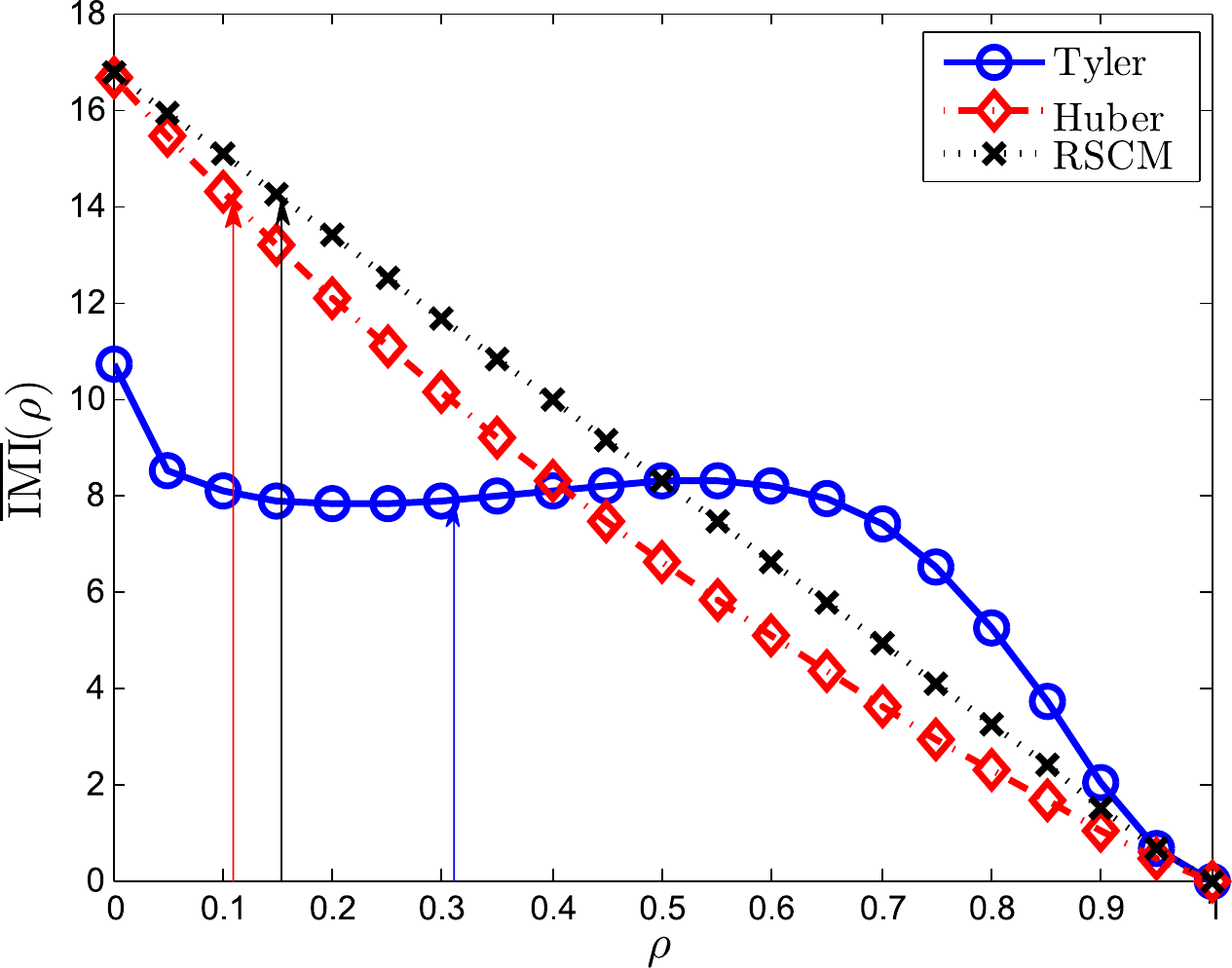}}
\end{minipage}
  \caption{{\small Infinitesimal measure of influence for $\rho \in (0,1]$ ($c_N = 3/2$). Arrows indicate the (oracle) optimal shrinkage parameter of the considered estimators in the absence of outliers. $[\mathbf{C}_N]_{ij}=.9^{|i-j|}$ and $[\mathbf{D}_N]_{ij}=.2^{|i-j|}$.}}
  \label{fig:if_reg1}
\end{figure}
Interestingly, the IMI of Tyler's estimator is somewhat less sensitive in $\rho$ than that of Huber's estimator or that of the RSCM (at least in the region where $\rho \approx \rho^\star$, indicated by arrows). This suggests that, in the presence of outliers, a small variation in the estimation of the optimal regularization parameter can have a different impact on the robustness of different estimators. Therefore, properly choosing the regularization parameter is crucial, both in terms of performance and robustness. These observations call for the need of a careful estimation of this optimal regularization parameter in the presence of outliers. 

Another important problem relates to the fact that, while the proposed (clean-data-optimal) choice of regularization parameter proves to be a practical
solution for a vanishing proportion of outliers, this choice would become more suboptimal under a more substantial non-vanishing
proportion of outliers. In such scenarios, and in particular if some a-priori knowledge on the outliers could be exploited,
different choices of regularization and/or estimators would be advisable. These problems will be investigated in future
work.

\appendix
\section{Proofs}
\subsection{Proofs of results in Section \ref{sec:shrinkage}}
\label{app:shrinkage}
\subsubsection{Theorem \ref{theorem1}}
\label{proof}
We will first start by proving the existence of $\hat{\mathbf{S}}_N(\rho)$, before turning to the continuity of $\rho \mapsto \gamma(\rho)$. Finally, we will show the uniform convergence of the spectral norm of $\hat{\mathbf{C}}_N(\rho)-\hat{\mathbf{S}}_N(\rho)$. The structure of the proof mirrors that of \cite[Theorem 1]{couillet2014large}, but non-trivial modifications are necessary to generalize the result {\color{black}to the 
class of $u$ functions} considered in this work. In particular, we have to make extensive use of the properties of the function $v=u(g^{-1}(x))$, where we recall that $g(x)=\frac{x}{1-(1-\rho)c\phi(x)}$.

Let us first prove the existence and uniqueness of $\hat{\mathbf{S}}_N(\rho)$. Notice that $\gamma$ (as defined in (\ref{eq:gamma_reg})) can be rewritten as the solution to the fixed-point equation
\begin{align}
\label{implicit_gamma}
\int \frac{t}{(1-\rho)\phi(g^{-1}(\gamma))t+\rho\gamma}\nu (dt) = 1.
\end{align}
Notice that the left-hand side of (\ref{implicit_gamma}) is a decreasing function of $\gamma$ (recall that $\phi$ is a increasing function and that $\rho>0$). Furthermore, it has limits $\infty$ as $\gamma \rightarrow 0$ (since $\phi(g^{-1}(0)) = 0$ and $\nu \ne \delta_0$ a.e.\@) and zero as $\gamma \rightarrow \infty$. This proves the existence and uniqueness of $\gamma$, from which the existence and uniqueness of $\hat{\mathbf{S}}_N(\rho)$ unfold.

Now, let us turn to the continuity of $\rho \mapsto \gamma(\rho)$. Consider a given compact set $I \in \mathcal{I}$, where $\mathcal{I}$ is the set of compact sets included in $ (0, 1]$. Take $\rho_1,\rho \in I$ and $\gamma_1=\gamma(\rho_1)$, $\gamma=\gamma(\rho)$. We can then write
\begin{align*}
&\int \frac{t}{(1-\rho)\phi(g^{-1}(\gamma))t+\rho\gamma}\nu (dt)  \\
&-\int \frac{t}{(1-\rho_1)\phi(g^{-1}(\gamma_1))t+\rho_1\gamma_1}\nu (dt) =0,
\end{align*}
which, after some algebra, leads to (\ref{eq:gamma_c}) (see top of next page).
\begin{figure*}[!]
\normalsize
\setcounter{equation}{21}
\begin{align}
\label{eq:gamma_c}
(\gamma_1-\gamma)\rho_1+\gamma(\rho_1-\rho)- \left((1-\rho)\phi(g^{-1}(\gamma))-(1-\rho_1)\phi(g^{-1}(\gamma_1))\right)\frac{\int\frac{t^2\nu(dt)}{((1-\rho)\phi(g^{-1}(\gamma))t+\rho\gamma)((1-\rho_1)\phi(g^{-1}(\gamma_1))t+\rho_1\gamma_1)}}{\int\frac{t\nu(dt)}{((1-\rho)\phi(g^{-1}(\gamma))t+\rho\gamma)((1-\rho_1)\phi(g^{-1}(\gamma_1))t+\rho_1\gamma_1)}}=0
\end{align}	
\setcounter{equation}{22}
\hrulefill
\end{figure*}
By Assumption 1.c, the support of $\nu$ is bounded by $\limsup_N ||\mathbf{C}_N|| < \infty$. In particular, recalling that $0 \le \phi(x) \le \phi_\infty$, from (\ref{implicit_gamma}) we necessarily have that $ \rho \gamma \le \limsup_N ||\mathbf{C}_N||$. It follows that the above integrals are uniformly bounded on $\rho$ in a neighborhood of $\rho_1 \le \rho$. Taking the limit $\rho \rightarrow \rho_1$, we then have (\ref{eq:gamma_c2}) (see top of next page).
\begin{figure*}[!]
\normalsize
\setcounter{equation}{22}
\begin{align}
\label{eq:gamma_c2}
(\gamma_1-\gamma)\rho_1+(1-\rho_1)\left(\phi(g^{-1}(\gamma_1))-\phi(g^{-1}(\gamma))\right)\frac{\int\frac{t^2\nu(dt)}{((1-\rho_1)\phi(g^{-1}(\gamma))t+\rho_1\gamma)((1-\rho_1)\phi(g^{-1}(\gamma_1))t+\rho_1\gamma_1)\nu(dt)}}{\int\frac{t\nu(dt)}{((1-\rho_1)\phi(g^{-1}(\gamma))t+\rho_1\gamma)((1-\rho_1)\phi(g^{-1}(\gamma_1))t+\rho_1\gamma_1)}} \rightarrow 0
\end{align}	
\setcounter{equation}{23}
\hrulefill
\end{figure*}
As $g^{-1}$ and $\phi$ are increasing, $\phi(g^{-1}(\gamma_1))-\phi(g^{-1}(\gamma))$ and $\gamma_1-\gamma$ have the same sign. As the above integrals are uniformly bounded on a neighborhood of $\rho_1$, we have $\gamma_1-\gamma \rightarrow 0$, from which we conclude that $\rho \mapsto \gamma(\rho)$ is continuous on $\mathcal{I}$.

Now, we will show the uniform convergence of the spectral norm of $\hat{\mathbf{C}}_N(\rho)-\hat{\mathbf{S}}_N(\rho)$. Let us fix $\rho \in I$, and denote $\hat{\mathbf{C}}_{(i)}(\rho)=\hat{\mathbf{C}}_N(\rho)-(1-\rho)\frac{1}{n}v\left(\frac{1}{N}\mathbf{y}_i^\dagger\hat{\mathbf{C}}_N^{-1}(\rho)\mathbf{y}_i\right)\mathbf{y}_i\mathbf{y}_i^\dagger$. After some algebra, we can rewrite
\[
\hat{\mathbf{C}}_N(\rho) = (1-\rho)\frac{1}{n}\sum_{i=1}^n v(d_i(\rho))\mathbf{y}_i\mathbf{y}_i^\dagger+\rho\mathbf{I}_N,
\]
where, for $i\in\{1,\cdots,n\}$, $d_i(\rho) \triangleq \frac{1}{N}\mathbf{y}_i^\dagger\hat{\mathbf{C}}_{(i)}^{-1}(\rho)\mathbf{y}_i$. Without loss of generality, we can assume that $d_1(\rho) \le \cdots \le d_n(\rho)$. Then, using the fact that $v$ is non-increasing, and the fact that $\mathbf{A} \succeq \mathbf{B} \Rightarrow \mathbf{B}^{-1} \succeq \mathbf{A}^{-1}$ for positive Hermitian matrices $\mathbf{A}$ and $\mathbf{B}$, we have
\begin{align*}
\begin{array}{rl}
d_n(\rho)&\hspace{-2.5mm}\!=\!\frac{1}{N}\mathbf{y}_n^\dagger\!\left(\!(1-\rho)\frac{1}{n}\sum_{j=1}^{n-1}v(d_j(\rho))\mathbf{y}_j\mathbf{y}_j^\dagger+\rho \mathbf{I}_N\!\right)^{-1}\!\mathbf{y}_n\\
&\hspace{-2.5mm} \!\le \!\frac{1}{N}\mathbf{y}_n^\dagger\!\left(\!(1-\rho)\frac{1}{n}\sum_{j=1}^{n-1}v(d_n(\rho))\mathbf{y}_j\mathbf{y}_j^\dagger+\rho \mathbf{I}_N\!\right)^{-1}\!\mathbf{y}_n.
\end{array}
\end{align*}
Since $\mathbf{y}_n \ne 0$ with probability $1$, we then have
\begin{align}
\label{dn}
& \mathbf{y}_n^\dagger \! \left(\!(1-\rho)\frac{1}{n}\sum_{j=1}^{n-1}d_n(\rho)v(d_n(\rho))\mathbf{y}_j\mathbf{y}_j^\dagger \! + \! \rho d_n(\rho) \mathbf{I}_N\!\right)^{-1}\!\!\!\mathbf{y}_n \nonumber\\
& \geq N.
\end{align}
Similarly,
\begin{align*}
&\mathbf{y}_1^\dagger\!\left(\!(1-\rho)\frac{1}{n}\sum_{j=2}^{n}d_1(\rho)v(d_1(\rho))\mathbf{y}_j\mathbf{y}_j^\dagger \! + \! \rho d_1(\rho) \mathbf{I}_N\!\right)^{-1}\!\mathbf{y}_1 \\
&\le N.
\end{align*}
We want to show that:
\begin{align*}
\sup_{\rho \in \mathcal{I}} \max_{1\le i\le n} \left|d_i(\rho)-\gamma(\rho)\right| \xrightarrow{\substack{\mathrm{a.s.}}} 0.
\end{align*}
This will be proven by a contradiction argument: assume there exists a sequence $\{\rho_n\}_{n=1}^\infty$ over which $d_n(\rho_n)>\gamma(\rho_n)+l$ infinitely often, for some $l>0$ fixed. Let us consider a subsequence of $\{\rho_n\}_{n=1}^\infty$ such that $\rho_n \rightarrow \rho_1$ (since $\{\rho_n\}_{n=1}^\infty$ is bounded, such subsequence exists by the Bolzano-Weierstrass theorem). On this subsequence, (\ref{dn}) gives us (\ref{en}) (see top of next page).
\begin{figure*}[!]
\normalsize
\setcounter{equation}{24}
\begin{align}
\label{en}
1  \le \frac{1}{N}\mathbf{y}_n^\dagger\left((1-\rho_n)\frac{1}{n}\sum_{j=1}^{n-1}(\gamma(\rho_n)+l)v(\gamma(\rho_n)+l)\mathbf{y}_j\mathbf{y}_j^\dagger+\rho_n (\gamma(\rho_n)+l)\mathbf{I}_N\right)^{-1}\mathbf{y}_n\triangleq e_n
\end{align}
\setcounter{equation}{25}
\hrulefill
\end{figure*}
Assume for now that $\rho_1 \ne 1$. Rewriting $xv(x)=\psi(x)$, we have (\ref{e+2}) (see top of next page), 
\begin{figure*}[!]
\normalsize
\setcounter{equation}{25}
\begin{align}
e_n & = \frac{1}{(1-\rho_n)\psi(\gamma(\rho_n)+l)}\frac{1}{N}\mathbf{y}_n^\dagger\left(\frac{1}{n}\sum_{j=1}^{n-1}\mathbf{y}_j\mathbf{y}_j^\dagger+\rho_n \frac{\gamma(\rho_n)+l}{(1-\rho_n)\psi(\gamma(\rho_n)+l)}\mathbf{I}_N\right)^{-1}\mathbf{y}_n \nonumber\\  
& \hspace{-2mm} \xrightarrow{\substack{\mathrm{a.s.}}} \frac{1}{(1-\rho_1)\psi(\gamma(\rho_1+l))} \delta\left(-(\gamma(\rho_1)+l)\rho_1\frac{1}{(1-\rho_1)\psi(\gamma(\rho_1)+l)}\right)\triangleq e^+
\label{e+2}
\end{align}
\setcounter{equation}{26}
\hrulefill
\end{figure*}
where, for $x>0$, $\delta(x)$ is the unique positive solution to the equation 
\begin{align}
\label{eq:delta}
\delta(x)= \int \frac{t}{-x+\frac{t}{1+c \delta(x)}}\nu(dt).
\end{align}
The convergence above follows from random matrix tools exposed in the proof of \cite[Theorem 1]{couillet2014large}. 
Define $(l,e) \mapsto h(l,e)$ as
\begin{align*}
h(l,e) \triangleq \int \frac{t}{(\gamma(\rho_1)+l)\rho_1 e+\frac{te}{\frac{1}{(1-\rho_1)\psi(\gamma(\rho_1)+l)}+ce}} \nu(dt),
\end{align*} 
which is clearly decreasing in both $l$ and $e$. Using (\ref{e+2}) and (\ref{eq:delta}) and a little algebra, we have that $h(l,e^+) =1$ for all $l >0$. Furthermore, from the definition of $\gamma(\rho_1)$, we also have that $h(0,1) = 1$. Therefore, $h(0,1)=h(l,e^+) =1$ for all $l >0$. Along with the fact that $e \mapsto h(\cdot,e)$ and $l \mapsto h(l,\cdot)$ are both decreasing functions, we then necessarily have $e^+ <1$.
But this is in contradiction with $e_n \ge 1$ from (\ref{en}). 

\noindent Assume now that $\rho_1=1$. Since $\frac{1}{N}||\mathbf{y}_n||^2 \xrightarrow{\substack{\mathrm{a.s.}}} M_{\nu,1} < \infty$, $\limsup_n ||\frac{1}{n}\sum_{i=1}^n \mathbf{y}_i\mathbf{y}_i^\dagger||< \infty$ a.s.\@ (from Assumption 1.b.\@ and \cite{bai1998no}), and $\gamma(1)=M_{\nu,1}$, from the definition of $e_n$ we have:
\[
e_n \xrightarrow{\substack{\mathrm{a.s.}}} \frac{M_{\nu,1}}{M_{\nu,1}+l}<1,
\]
which is again a contradiction.

It follows that for all large $n$, there is no sequence of $\rho_n$ such that $d_n(\rho)> \gamma(\rho)+l$ infinitely often. Consequently, $d_n(\rho) \le \gamma(\rho)+l$ for all large $n$ a.s., uniformly on $\rho \in I$. 
We can apply the same strategy to prove that $d_1(\rho)$ is greater than $\gamma(\rho)-l$ for all large $n$ uniformly on $\rho \in I$. As this is true for arbitrary $l>0$, we then have $\sup_{\rho \in I} \max_{1 \le i \le n} |d_i(\rho)-\gamma(\rho)|\xrightarrow{\substack{\mathrm{a.s.}}} 0$. By continuity of $v$, we also have  $\sup_{\rho \in I} \max_{1 \le i \le n} |v(d_i(\rho))-v(\gamma(\rho))|\xrightarrow{\substack{\mathrm{a.s.}}} 0$.
It follows that
\vspace{5mm}
\begin{align*}
&\sup_{\rho \in I} \left|\left|\hat{\mathbf{C}}_N(\rho)-\hat{\mathbf{S}}_N(\rho)\right|\right| \\
& \le \left|\left|\frac{1}{n}\sum_{i=1}^n\mathbf{y}_i\mathbf{y}_i^\dagger\right|\right| \sup_{\rho \in I} \max_{1\le i \le n}(1-\rho) |v(d_i)-v(\gamma)| \xrightarrow{\substack{\mathrm{a.s.}}} 0,
\end{align*}
where we used the fact that $\limsup_n \left|\left|\frac{1}{n}\sum_{i=1}^n \mathbf{y}_i\mathbf{y}_i^\dagger \right|\right| < \infty$ a.s., as above.

\subsubsection{Proposition \ref{cor 2}}
Since $\rho \mapsto v(\gamma)$ is non-negative, it is clear that $\underline{\rho}$ is indeed in $(0,1]$.  Then, for a couple $(\rho,\underline{\rho})$ satisfying (\ref{eq:rho_rho_bar}), the (relative) weights given to the SCM $\frac{1}{n}\sum_{i=1}^n \mathbf{y}_i \mathbf{y}_i^\dagger$ and the shrinkage target $\mathbf{I}_N$ are the same for $\hat{\mathbf{S}}_N(\rho)$ and $\mathbf{R}(\underline{\rho})$. After trace-normalization, the first result of Proposition \ref{cor 2} follows. Now, since $\rho \mapsto v(\gamma)$ is continuous and bounded (from Theorem \ref{theorem1}), it follows that $F: (0,1] \rightarrow (0,1]$ is continuous and onto, from which the second result of Proposition \ref{cor 2} unfolds. 

\subsubsection{Proposition \ref{cor 3}}
The proof of Proposition \ref{cor 3} makes use of the asymptotic equivalence of $\hat{\mathbf{C}}_N(\rho)$ with $\hat{\mathbf{S}}_N(\rho)$ (as given in Theorem \ref{theorem1}) and the equivalence and mapping between $\hat{\mathbf{S}}_N(\rho)$ and the RSCM (as given in Proposition \ref{cor 2}). It is known that the RSCM can be optimized with respect to the Frobenius norm, with the corresponding optimal regularization parameter $\rho^\star$ given in Proposition \ref{cor 3} (see, e.g., \cite{ledoit2004well,couillet2014large}). With (\ref{eq:rho_rho_bar}), it follows that for $\hat{\rho}^\star$ a solution to $\frac{\hat{\rho}^\star}{(1-\hat{\rho}^\star)  v(\gamma)+\hat{\rho}^\star}=\rho^\star$, the associated estimator $\hat{\mathbf{C}}_N(\hat{\rho}^\star)$ will have (asymptotically) minimal quadratic loss. Similarly to \cite[Proposition 2]{couillet2014large}, the second part of Proposition \ref{cor 3} provides a consistent estimate $\hat{\rho}^\star$ based on a possible estimate of $\rho^\star$, the optimal regularization parameter for the RSCM.

\vspace{-5mm}

\subsection{Proofs of the results in Section \ref{sec:if}}
\label{app:if1}
\subsubsection{Theorem \ref{theorem_no_shrinkage}}  
The convergence of the spectral norm of $\hat{\mathbf{C}}_N^{\epsilon_n}-\hat{\mathbf{S}}_N^{\epsilon_n}$ unfolds from the proof of \cite[Theorem 1]{morales2015large}. However, the proof of the existence and uniqueness of $\hat{\mathbf{S}}_N^{\epsilon_n}$ for $\epsilon$ arbitrary requires additional arguments. To proceed, we make use of the standard interference function framework \cite{yates1995framework}. Define the real-valued functions $h_i:[0,\infty) \rightarrow [0,\infty),(q_0,q_1) \mapsto h_i(q_0,q_1)$ (with $i=0,1$) as:
\begin{align*}
h_0(q_0,q_1) & = \frac{1}{N} \! \tr  \mathbf{C}_N \! \left(\!(1-\epsilon) \frac{f(q_0)}{q_0}\mathbf{C}_N\! + \!\epsilon \frac{f(q_1)}{q_1} \mathbf{D}_N\!\right)^{-1} \\
h_1(q_0,q_1) & = \frac{1}{N} \! \tr  \mathbf{D}_N \! \left(\!(1-\epsilon) \frac{f(q_0)}{q_0}\mathbf{C}_N\! +\! \epsilon \frac{f(q_1)}{q_1} \mathbf{D}_N \! \right)^{-1} \!,
\end{align*}
where $f(x) \triangleq  \frac{xv(x)}{1+c xv(x)}$, and where we dropped the subscript $n$ of $\epsilon_n$ for readability. Thus defined, $f$ is onto from $[0,\infty)$ to $[0,\phi_\infty)$, where we recall that $\phi_\infty > 1$. 
It can be easily verified that $h_0$, $h_1$ are standard interference functions\footnote{In particular, they should verify conditions of positivity, monotonocity and scalability.} (see \cite{morales2015large} for details). According to \cite[Theorem 2]{yates1995framework}, if there exist some $q_0, q_1 >0$ such that $h_0(q_0,q_1) \le q_0$ and $h_1(q_0,q_1) \le q_1$, then the system of fixed-point equations $h_0(q_0,q_1) = q_0$, $h_1(q_0,q_1) = q_1$ admits a unique solution $\{q_0,q_1\}$. It therefore remains to find $q_0$ and $q_1$ that satisfy $h_0(q_0,q_1) \le q_0$ and $h_1(q_0,q_1) \le q_1$. {\color{black} If we were to assume, like in \cite{morales2015large}, that $\epsilon <1-c$, with also $(1-\epsilon)^{-1} < \phi_\infty < c^{-1}$, then we could find explicit $q_0, q_1$ depending on $\mathbf{C}_N, \mathbf{D}_N$, and the function $v$ (see \cite{morales2015large} for details). However, in the more general case where $\epsilon \in [0,1)$, such $q_0, q_1$ can not be found explicitly, making the derivation more challenging. To prove the existence of eligible $q_0, q_1$, we have to make use of the relationship between the equations $h_0$ and $h_1$. It turns out that it is sufficient to show the existence of $q_0,q_1\ge f^{-1}(1)$ such that }
\begin{align*}
h'_0(q_0,q_1) &\le q_0 \\
h'_1(q_0,q_1) &\le q_1,
\end{align*}
where
\begin{align*}
h'_0(q_0,q_1) &\triangleq  \frac{1}{N} \tr \mathbf{C}_N\left((1-\epsilon) \frac{1}{q_0}\mathbf{C}_N + \epsilon \frac{1}{q_1} \mathbf{D}_N\right)^{-1} \\
& \geq h_0(q_0,q_1)\\
h'_1(q_0,q_1) &\triangleq  \frac{1}{N} \tr \mathbf{D}_N\left((1-\epsilon) \frac{1}{q_0}\mathbf{C}_N + \epsilon \frac{1}{q_1} \mathbf{D}_N\right)^{-1}  \\
& \geq h_1(q_0,q_1).
\end{align*}

Consider two cases depending on whether $\epsilon \in \{0,1\}$ or not.
If $\epsilon = 0$, then 
\begin{align*}
h'_0(q_0,q_1) & = q_0 \\
h'_1(q_0,q_1) & = \frac{1}{N} \tr \mathbf{D}_N\mathbf{C}_N^{-1} q_0.
\end{align*}
Taking $q_1 = \mathbf{D}_N\mathbf{C}_N^{-1} q_0$, we have $h'_0(q_0,q_1) \le q_i$ for $i=0,1$. It remains to choose $q_0$ such that $\min \{q_0, q_1\} \ge f^{-1}(1)$ (which is always possible), and the proof is done. 
Similarly, if $\epsilon = 1$, it suffices to take $q_0 = \mathbf{C}_N\mathbf{D}_N^{-1} q_1$, with $q_1$ chosen such that $\min \{q_0, q_1\} \ge f^{-1}(1)$.

Assume now that $\epsilon \in (0,1)$.
Let us define $\alpha \triangleq \frac{q_1}{q_0} > 0$. We can rewrite:
\begin{align*}
h'_0(q_0,q_1) &= q_0  \frac{1}{N} \tr \left(\!(1-\epsilon) \mathbf{I}_N +  \frac{\epsilon}{\alpha} \mathbf{C}_N^{-1}\mathbf{D}_N\!\right)^{-1}\\
h'_1(q_0,q_1) &= \frac{q_1}{\alpha} \frac{1}{N} \tr \mathbf{C}_N^{-1}\mathbf{D}_N\!\left(\!(1-\epsilon) \mathbf{I}_N +  \frac{\epsilon}{\alpha} \mathbf{C}_N^{-1}\mathbf{D}_N \!\right)^{-1} .
\end{align*}
Finding $q_0, q_1$ such that $h'_i \le q_i$ is then equivalent to finding $\alpha$ such that 
\begin{align}
\frac{1}{N} \tr \left((1-\epsilon) \mathbf{I}_N + \epsilon \frac{1}{\alpha} \mathbf{C}_N^{-1}\mathbf{D}_N\right)^{-1} \le 1  \label{eq:h11}\\
\frac{1}{\alpha}\frac{1}{N} \tr \mathbf{C}_N^{-1}\mathbf{D}_N \left((1-\epsilon) \mathbf{I}_N + \epsilon \frac{1}{\alpha} \mathbf{C}_N^{-1}\mathbf{D}_N\right)^{-1} \le 1. \label{eq:h12}
\end{align}
By applying Lemma \ref{lem:equiv} (see below) with $\mathbf{A}=\frac{1}{\alpha}\mathbf{C}_N^{-1}\mathbf{D}_N$, we can show that
\begin{align*}
\frac{1}{\alpha}\frac{1}{N} \tr \mathbf{C}_N^{-1}\mathbf{D}_N \left((1-\epsilon) \mathbf{I}_N + \epsilon \frac{1}{\alpha} \mathbf{C}_N^{-1}\mathbf{D}_N\right)^{-1} &\le 1 \\
\Leftrightarrow \frac{1}{N} \tr \left((1-\epsilon) \mathbf{I}_N + \epsilon \frac{1}{\alpha} \mathbf{C}_N^{-1}\mathbf{D}_N\right)^{-1} &\ge 1.
\end{align*}
Combined with (\ref{eq:h11}) and (\ref{eq:h12}), it follows that $q_0, q_1$ verify $h'_i \le q_i$ if and only if 
\begin{align*}
\frac{1}{N} \tr \left((1-\epsilon) \mathbf{I}_N + \epsilon \frac{1}{\alpha} \mathbf{C}_N^{-1}\mathbf{D}_N\right)^{-1} = 1.
\end{align*}

Denote by $a_i>0$ the $i$-th eigenvalue of $\mathbf{C}^{-1}\mathbf{D}$. We then have:
\begin{align*}
&\frac{1}{N} \tr \left((1-\epsilon) \mathbf{I}_N + \epsilon \frac{1}{\alpha} \mathbf{C}_N^{-1}\mathbf{D}_N\right)^{-1} = 1 \\
& \Leftrightarrow \frac{1}{N}\sum_{i=1}^N \frac{1}{1-\epsilon+\epsilon \frac{a_i}{\alpha}} = 1 \\
& \Leftrightarrow \frac{1}{N}\alpha\sum_{i=1}^N \prod_{j \neq i} \left((1-\epsilon)\alpha+\epsilon a_j\right) = \prod_{i=1}^N\left((1-\epsilon)\alpha +\epsilon a_i\right),
\end{align*}
where the last equality comes from putting all the terms in the sum on the same denominator, and multiplying by $\alpha^N \neq 0$. Finding an eligible $\alpha$ therefore boils down to finding whether the polynomial in $\alpha$ appearing in the last equation has positive roots. Notice now that the leading coefficient of this $N$-order polynomial is $b_N = \epsilon(1-\epsilon)^{N-1}>0$,
while the constant is $b_0 = -\epsilon^N \prod_{i=1}^N a_i < 0$.
As $b_N \times b_0 <0$, it follows that this polynomial admits (at least) one positive root (by applying the intermediate value theorem). Call $\alpha_0$ such a root. Choosing a $q_0$ that verifies $\min \{q_0, q_0\alpha_0 \} \ge f^{-1}(1)$, $q_0$ and taking $q_1 =q_0\alpha_0$, we will then have $h'_i \le q_i$, and therefore $h_i \le q_i$. The existence and uniqueness of $\gamma^\epsilon$ and $\alpha$, as given in Theorem~\ref{theorem_no_shrinkage}, unfold.

\begin{lemma} For $\mathbf{A}$ an invertible matrix and $\alpha > 0$, $\epsilon \in (0,1)$, we have the following equivalence:
\label{lem:equiv} 
\begin{align*}
\frac{1}{N} \tr \mathbf{A}\left((1-\epsilon) \mathbf{I}_N + \epsilon \mathbf{A}\right)^{-1} \le 1 &\Leftrightarrow \\
\frac{1}{N} \tr \left((1-\epsilon) \mathbf{I}_N + \epsilon \mathbf{A}\right)^{-1} \ge 1&.
\end{align*}
\end{lemma}

\subsubsection{Corollary \ref{cor:mi_nonreg}}
This is a direct consequence of Theorems~\ref{th:Romain_no_shrinkage} and \ref{theorem_no_shrinkage}, by writing $\overline{\mathrm{MI}}(\epsilon_n) =\left\lVert \mathbb{E} \left[\frac{\hat{\mathbf{S}}_N^{0}}{L^0}-\frac{\hat{\mathbf{S}}_N^{\epsilon_n}}{L^\epsilon}\right] \right\rVert$, where $L^{\epsilon_n}_N \triangleq v(\gamma^{\epsilon_n})(1-\epsilon_n)  + v(\alpha^{\epsilon_n})\epsilon_n $, and using the fact that $\frac{1}{N}\tr \hat{\mathbf{C}}_N^{\epsilon_n} - L^{\epsilon_n}_N \xrightarrow{\substack{\mathrm{a.s.}}} 0$ and $\frac{1}{N}\tr \hat{\mathbf{S}}_N^{\epsilon_n} - L^{\epsilon_n}_N \xrightarrow{\substack{\mathrm{a.s.}}} 0$. 

\subsection{Proof of the results in Section \ref{sec:if_reg}}
\label{app:if2}
\subsubsection{Theorem \ref{theorem2}}
The convergence of the spectral norm of $\hat{\mathbf{C}}_N^{\epsilon_n}(\rho)-\hat{\mathbf{S}}_N^{\epsilon_n}(\rho)$ for $\rho \in \mathcal{R}$ is a direct extension of the proof of \cite[Theorem 1]{morales2015large}, adapted to account for the introduction of the regularization parameter $\rho \in \mathcal{R}$. 
{\color{black} We now turn to the proof of existence and uniqueness of $\hat{\mathbf{S}}_N^{\epsilon_n}(\rho)$. To proceed, we want to prove that the system of equations (\ref{eq:rd_shrinkage1}) in Theorem~\ref{theorem2} admits a unique solution $\{\gamma, \alpha\}$ for a fixed $\rho \in \mathcal{R}$ (with $\mathcal{R}$ defined at the beginning of Section \ref{sec:out}). As for the proof of Theorem~\ref{theorem_no_shrinkage}, we make use of the standard interference function framework \cite{yates1995framework}. Define the real-valued functions $h^\rho_i:[0,\infty) \rightarrow [0,\infty),(q_0,q_1) \mapsto h^\rho_i(q_0,q_1)$
\begin{align}
\label{eq:h0}
h^\rho_0(q_0,q_1) & = \frac{1}{N}  \tr  \mathbf{C}_N  \mathbf{E}_N^{-1} \\
\label{eq:h1}
h^\rho_1(q_0,q_1) & = \frac{1}{N}  \tr  \mathbf{D}_N  \mathbf{E}_N^{-1},
\end{align}
where
\begin{align*}
\mathbf{E}_N \!= \!\left(\!\!(1\!-\!\rho)(1\!-\!\epsilon) \frac{f^\rho(q_0)}{q_0}\mathbf{C}_N \!+\!(1\!-\!\rho)\epsilon \frac{f^\rho(q_1)}{q_1} \mathbf{D}_N \!+\!\rho \mathbf{I}_N\!\!\right),
\end{align*}
with $f^\rho(x) \triangleq  \frac{xv(x)}{1+c (1-\rho) xv(x)}$. Following the same argument as for the proof of Theorem~\ref{theorem_no_shrinkage}, we want to find some $q_0$, $q_1$ such that $h^\rho_i(q_0,q_1) \le q_i$ (with $i=0,1$) for all $\rho \in \mathcal{R}$. 
Fortunately, in contrast to the non-regularized case, we can find such $q_0, q_1$ explicitly, thanks to the regularization parameter $\rho$. 
Indeed, since $\mathbf{E}_N\succeq  \rho \mathbf{I}_N,$ it follows that, from \cite[Corollary 7.7.4]{horn2012matrix}, $\mathbf{E}_N^{-1} \preceq \frac{1}{\rho}\mathbf{I}_N$,
with $\rho >0$.
Taking $q_0 = \frac{1}{\rho}\frac{1}{N}\mathrm{Tr} \mathbf{C}_N$ and $q_1 = \frac{1}{\rho}\frac{1}{N}\mathrm{Tr} \mathbf{D}_N$, and plugging in (\ref{eq:h0}), (\ref{eq:h1}), we have
\begin{align*}
h^\rho_0(q_0,q_1) &\le q_0 \\
h^\rho_1(q_0,q_1) &\le q_1,
\end{align*}
which, using \cite[Theorem 2]{yates1995framework}, proves that there exist $q_0, q_1$ such that $h^\rho_i(q_0,q_1) = q_i$ for $i=0,1$. This proves existence and uniqueness of $\hat{\mathbf{S}}_N^{\epsilon_n}(\rho)$. 
}

\subsubsection{Corollary \ref{cor:mi_shrinkage}}
Define 
\begin{align*}
\overline{\mathrm{MI}}(\rho,\epsilon_n) \triangleq \left\lVert \mathbb{E}\left[ \frac{\hat{\mathbf{S}}_{N}^{0}(\rho)}{L^0(\rho)}-\frac{\hat{\mathbf{S}}_{N}^{\epsilon_n}(\rho)}{L^{\epsilon_n}(\rho)}\right] \right\rVert,
\end{align*}
with $L^{\epsilon_n}_N(\rho) \triangleq v(\gamma^{\epsilon_n})(1-\rho)(1-\epsilon_n)  + v(\alpha^{\epsilon_n})(1-\rho)\epsilon_n  + \rho$. The convergence result is a direct consequence of Theorem~\ref{theorem2} and the fact that, for $\rho \in \mathcal{R}$, $\frac{1}{N}\tr \hat{\mathbf{C}}_N^{\epsilon_n}(\rho) - L^{\epsilon_n}_N(\rho) \xrightarrow{\substack{\mathrm{a.s.}}} 0$ and $\frac{1}{N}\tr \hat{\mathbf{S}}_N^{\epsilon_n}(\rho)- L^{\epsilon_n}_N(\rho) \xrightarrow{\substack{\mathrm{a.s.}}} 0$. The derivation of $\overline{\mathrm{MI}}(\rho,\epsilon_n) $ is straightforward by expanding
\begin{align*}
\overline{\mathbf{L}}(\rho,\epsilon_n) & \triangleq \mathbb{E}\left[ \frac{\hat{\mathbf{S}}_{N}^{0}(\rho)}{L^0(\rho)}-\frac{\hat{\mathbf{S}}_{N}^{\epsilon_n}(\rho)}{L^{\epsilon_n}(\rho)}\right]  \\
& \hspace{-13mm} =  \frac{(1-\rho)(1-\epsilon_n)v(\gamma^{\epsilon_n})\mathbf{C}_N+(1-\rho)\epsilon_n v(\alpha^{\epsilon_n})\mathbf{D}_N+\rho \mathbf{I}_N}{(1-\rho)(1-\epsilon_n)v(\gamma^{\epsilon_n})+(1-\rho)\epsilon_n v(\alpha^{\epsilon_n}) + \rho} \\
& \hspace{-9mm}-\frac{(1-\rho)v(\gamma^0)\mathbf{C}_N+\rho \mathbf{I}_N}{(1-\rho)v(\gamma^0)+\rho} \\
& \hspace{-13mm} =\frac{U(\epsilon_n,\rho)}{V(\epsilon_n,\rho)}, 
\end{align*}
with $U(\epsilon_n,\rho)$ and $V(\epsilon_n,\rho)$ given in the corollary.

\subsubsection{Corollary \ref{cor:ia_shrinkage}}
The result follows directly by taking the limit of $\frac{U(\epsilon_n,\rho)}{V(\epsilon_n,\rho)}$, as given in Corollary \ref{cor:mi_shrinkage}. 

\subsubsection{Computation of $\overline{\mathrm{IMI}}(\rho)$} 
\label{app:if}
In order to compute $\overline{\mathrm{IMI}}(\rho)$ (\ref{eq:if_shrinkage}) for arbitrary $\rho$, we need to compute $\big. \frac{dv}{d\epsilon_n} \big|_{\epsilon_n = 0}=\big. \frac{dv}{d\gamma} \big|_{\gamma = \gamma^0} \times \big. \frac{d\gamma}{d\epsilon_n} \big|_{\epsilon_n = 0}$. For this, let us adopt the following notations:
\begin{align*}
\mathbf{A}_N(\rho) &= \frac{(1-\rho)v(\gamma^0)}{1+(1-\rho)c\gamma^0v(\gamma^0)}\mathbf{C}_N+\rho \mathbf{I}_N \\
\gamma^0 & \triangleq \lim_{n \rightarrow \infty} \gamma^{\epsilon_n} = \frac{1}{N} \tr \mathbf{C}_N \mathbf{A}_N^{-1}(\rho) \\
\alpha^0 & \triangleq \lim_{n \rightarrow \infty} \alpha^{\epsilon_n} = \frac{1}{N} \tr \mathbf{D}_N \mathbf{A}_N^{-1}(\rho).
\end{align*}

Let us first compute $\big. \frac{d\gamma}{d\epsilon_n} \big|_{\epsilon_n = 0}$. For this, we need to differentiate (\ref{eq:rd_shrinkage1}) with respect to $\epsilon_n$. We can do so by using the fact that $\frac{d\mathbf{M}^{-1}}{d\zeta}(\zeta)=-\mathbf{M}^{-1}(\zeta)\frac{d\mathbf{M}}{d\zeta}(\zeta)\mathbf{M}^{-1}(\zeta)$. Taking the limit when $\epsilon_n \rightarrow 0$ in the resulting equation, we get (\ref{eq:dvde}) (see top of next page).
\begin{figure*}[!]
\normalsize
\setcounter{equation}{31}
\begin{align}
\label{eq:dvde}
\left. \frac{d\gamma}{d\epsilon_n}\right|_{\epsilon_n = 0} = (1-\rho)\frac{\frac{1}{N} \tr \left[ \mathbf{A}_N^{-1}(\rho) \mathbf{C}_N\mathbf{A}_N^{-1}(\rho)\left(\frac{v(\gamma^0)}{1+(1-\rho)c \gamma^0 v(\gamma^0)}\mathbf{C}_N -\frac{v(\alpha^0)}{1+(1-\rho)c \alpha^0 v(\alpha^0)} \mathbf{D}_N\right) \right] }{1+\frac{(1-\rho)\left(\left. \frac{dv}{d\gamma}\right|_{\epsilon = 0}-(1-\rho)c v(\gamma^0)^2\right)}{(1+(1-\rho)c\gamma^0 v(\gamma^0))^2} \frac{1}{N} \tr \mathbf{A}_N^{-1}(\rho) \mathbf{C}_N\mathbf{A}_N^{-1}(\rho)\mathbf{C}_N}
\end{align}	
\setcounter{equation}{32}
\hrulefill
\end{figure*}
It remains to compute $\big. \frac{dv}{d\gamma} \big|_{\gamma = \gamma^0}$. It is challenging to find a general expression for $\big. \frac{dv}{d\gamma} \big|_{\gamma = \gamma^0}$ for an arbitrary $u$ function (since it requires computing $v(x) = u(g^{-1}(x))$, which does not necessarily take a tractable form). However, we can do so for the $u$ functions $u_\mathrm{M-Tyler}(x)=\frac{1}{c_N}\frac{1+t}{t+x}$ and $u_\mathrm{M-Huber}(x)=\frac{1}{c_N}\min \{1,\frac{1+t}{t+x} \}$. Assume $\rho > 0$ (for $\rho =0$ (when possible), we fall back into the non-regularized case). 
Then, for these two functions, the associated $v$ functions can be approximated by
 \begin{align*}
v_\mathrm{M-Tyler}(x) & \simeq \frac{1}{c_N}\frac{1+t}{t+\rho x} 
\end{align*}
and
 \begin{align*}
v_\mathrm{M-Huber}(x) & \simeq \left\{ \begin{array}{ll}
  \frac{1}{c_N} & \text{\quad if \quad} x \leq \frac{1}{\rho} \\
  \frac{1}{c_N}\frac{1+t}{t+\rho x}  & \text{\quad if \quad} x \geq \frac{1}{\rho}
 \end{array}  \right.,
\end{align*}
for $t$ small, from which we can deduce $\big. \frac{dv}{dx} \big|_{x = \gamma^0}$.\footnote{Note however that $v_\mathrm{M-Huber}$ is only piece-wise differentiable. In particular, additional care is needed if $\gamma^0 = 1/\rho$.} We can then substitute $\big. \frac{dv}{d\epsilon}\big|_{\epsilon = 0}= \big. \frac{dv}{d\gamma}\big|_{\gamma = \gamma^0} \times \big. \frac{d\gamma}{d\epsilon}\big|_{\epsilon = 0}$ in (\ref{eq:if_shrinkage}). We can proceed similarly for $u_\mathrm{M-Tyler}(x)=\frac{1+t}{t+x}$ and $u_\mathrm{M-Huber}(x)=\min \{1,\frac{1+t}{t+x} \}$.

\balance
\bibliographystyle{ieeetr}
\small
\bibliography{References}

\begin{thebibliography}{10}

\bibitem{abramovich2007diagonally}
Y.~Abramovich and N.~K. Spencer, ``Diagonally loaded normalised sample matrix
  inversion ({LNSMI}) for outlier-resistant adaptive filtering,'' in {\em IEEE
  Int. Conf. Acoust. Signal Process.}, vol.~3, pp.~III--1105, 2007.

\bibitem{pascal2014generalized}
F.~Pascal, Y.~Chitour, and Y.~Quek, ``Generalized robust shrinkage estimator
  and its application to {STAP} detection problem,'' {\em IEEE Trans. Signal
  Process.}, vol.~62, pp.~5640--5651, Sept. 2014.

\bibitem{tulino2004random}
A.~M. Tulino and S.~Verd{\'u}, ``Random matrix theory and wireless
  communications,'' {\em Foundations and Trends in Communications and
  Information Theory}, vol.~1, no.~1, pp.~1--182, 2004.

\bibitem{ledoit2003improved}
O.~Ledoit and M.~Wolf, ``Improved estimation of the covariance matrix of stock
  returns with an application to portfolio selection,'' {\em J. Empir.
  Finance}, vol.~10, no.~5, pp.~603--621, 2003.

\bibitem{schafer2005shrinkage}
J.~Sch{\"a}fer and K.~Strimmer, ``A shrinkage approach to large-scale
  covariance matrix estimation and implications for functional genomics,'' {\em
  Stat. {A}pplicat. {G}enetics {M}olecular {B}iology}, vol.~4, no.~1, 2005.

\bibitem{ward1981compound}
K.~D. Ward, ``Compound representation of high resolution sea clutter,'' {\em
  Electron. Lett.}, vol.~17, no.~16, pp.~561--563, 1981.

\bibitem{billingsley1999statistical}
J.~B. Billingsley, A.~Farina, F.~Gini, M.~V. Greco, and L.~Verrazzani,
  ``Statistical analyses of measured radar ground clutter data,'' {\em IEEE
  Trans. Aerosp. Electron. Syst.}, vol.~35, pp.~579--593, Apr. 1999.

\bibitem{kelker1970distribution}
D.~Kelker, ``Distribution theory of spherical distributions and a
  location-scale parameter generalization,'' {\em Sankhy{\=a}: Indian J. Stat.,
  Series A}, pp.~419--430, 1970.

\bibitem{ollila2012complex}
E.~Ollila, D.~E. Tyler, V.~Koivunen, and H.~V. Poor, ``Complex elliptically
  symmetric distributions: Survey, new results and applications,'' {\em IEEE
  Trans. Signal Process.}, vol.~60, pp.~5597--5625, Aug. 2012.

\bibitem{mestre2008modified}
X.~Mestre and M.~{\'A}. Lagunas, ``Modified subspace algorithms for {D}o{A}
  estimation with large arrays,'' {\em IEEE Trans. Signal Process.}, vol.~56,
  pp.~598--614, Jan. 2008.

\bibitem{nadler2010nonparametric}
B.~Nadler, ``Nonparametric detection of signals by information theoretic
  criteria: performance analysis and an improved estimator,'' {\em IEEE Trans.
  Signal Process.}, vol.~58, pp.~2746--2756, Feb. 2010.

\bibitem{huber1964robust}
P.~J. Huber, ``Robust estimation of a location parameter,'' {\em Ann. Math.
  Stat.}, vol.~35, no.~1, pp.~73--101, 1964.

\bibitem{maronna1976robust}
R.~A. Maronna, ``Robust {M}-estimators of multivariate location and scatter,''
  {\em {A}nn. {S}tat.}, pp.~51--67, 1976.

\bibitem{tyler1987distribution}
D.~E. Tyler, ``A distribution-free {M}-estimator of multivariate scatter,''
  {\em Ann. Stat.}, pp.~234--251, 1987.

\bibitem{chen2011robust}
Y.~Chen, A.~Wiesel, and A.~O. Hero, ``Robust shrinkage estimation of
  high-dimensional covariance matrices,'' {\em IEEE Trans. Signal Process.},
  vol.~59, pp.~4097--4107, Apr. 2011.

\bibitem{couillet2013robust}
R.~Couillet, F.~Pascal, and J.~W. Silverstein, ``Robust estimates of covariance
  matrices in the large dimensional regime,'' {\em IEEE Trans. Inf. Theory},
  vol.~60, pp.~7269--7278, Sept. 2014.

\bibitem{couillet2014large}
R.~Couillet and M.~R. McKay, ``Large dimensional analysis and optimization of
  robust shrinkage covariance matrix estimators,'' {\em J. {M}ultivar.
  {A}nal.}, vol.~131, pp.~99--120, Oct. 2014.

\bibitem{zhang2014marchenko}
T.~Zhang, X.~Cheng, and A.~Singer, ``Marchenko-{P}astur law for {T}yler's and
  {M}aronna's {M}-estimators,'' {\em arXiv preprint arXiv:1401.3424}, 2014.

\bibitem{morales2015large}
D.~Morales-Jimenez, R.~Couillet, and M.~R. McKay, ``Large dimensional analysis
  of robust {M}-estimators of covariance with outliers,'' {\em IEEE Trans.
  Signal Process.}, vol.~63, pp.~5784--5797, Jul. 2015.

\bibitem{abramovich1981controlled}
Y.~Abramovich, ``A controlled method for adaptive optimization of filters using
  the criterion of maximum signal-to-noise ratio,'' {\em Radio Eng. Elect.
  Phys}, vol.~26, no.~3, pp.~87--95, 1981.

\bibitem{carlson1988covariance}
B.~D. Carlson, ``Covariance matrix estimation errors and diagonal loading in
  adaptive arrays,'' {\em IEEE Trans. Aerosp. Electron. Syst.}, vol.~24,
  pp.~397--401, Jul. 1988.

\bibitem{ollila2014regularized}
E.~Ollila and D.~E. Tyler, ``Regularized {M}-estimators of scatter matrix,''
  {\em IEEE Trans. Signal Process.}, vol.~62, pp.~6059--6070, Sept. 2014.

\bibitem{couillet2015random}
R.~Couillet, F.~Pascal, and J.~W. Silverstein, ``The random matrix regime of
  {M}aronna's {M}-estimator with elliptically distributed samples,'' {\em J.
  {M}ultivar. {A}nal.}, vol.~139, pp.~56--78, Jul. 2015.

\bibitem{kent1991redescending}
J.~T. Kent and D.~E. Tyler, ``Redescending {M}-estimates of multivariate
  location and scatter,'' {\em Ann. Stat.}, pp.~2102--2119, 1991.

\bibitem{soloveychik2014tyler}
I.~\color{black} Soloveychik and A.~Wiesel, ``Tyler's covariance matrix
  estimator in elliptical models with convex structure,'' {\em IEEE Trans.
  Signal Process.}, vol.~62, no.~20, pp.~5251--5259, 2014\color{black}.

\bibitem{yang2015robust}
L.~\color{black} Yang, R.~Couillet, and M.~R. McKay, ``A robust statistics
  approach to minimum variance portfolio optimization,'' {\em IEEE Trans.
  Signal Process.}, vol.~63, no.~24, pp.~6684--6697, 2015\color{black}.

\bibitem{couillet2016second}
R.~\color{black}Couillet, A.~Kammoun, and F.~Pascal, ``Second order statistics
  of robust estimators of scatter. application to {G}{L}{R}{T} detection for
  elliptical signals,'' {\em J. Multivar. Anal.}, vol.~143, pp.~249--274,
  2016\color{black}.

\bibitem{sun2014regularized}
Y.~Sun, P.~Babu, and D.~P. Palomar, ``Regularized {T}yler's scatter estimator:
  Existence, uniqueness, and algorithms,'' {\em IEEE Trans. Signal Process.},
  vol.~62, pp.~5143--5156, Aug. 2014.

\bibitem{huber2011robust}
P.~J. Huber, {\em Robust {S}tatistics}.
\newblock Springer, 2011.

\bibitem{ledoit2004well}
O.~Ledoit and M.~Wolf, ``A well-conditioned estimator for large-dimensional
  covariance matrices,'' {\em J. {M}ultivar. {A}nal.}, vol.~88, pp.~365--411,
  Jul. 2004.

\bibitem{du2010fully}
L.~Du, J.~Li, and P.~Stoica, ``Fully automatic computation of diagonal loading
  levels for robust adaptive beamforming,'' {\em IEEE Trans. Aerosp. Electron.
  Syst.}, vol.~46, Feb. 2010.

\bibitem{hoerl1970ridge}
A.~E. Hoerl and R.~W. Kennard, ``Ridge regression: {B}iased estimation for
  nonorthogonal problems,'' {\em Technometrics}, vol.~12, no.~1, pp.~55--67,
  1970.

\bibitem{kaufman2009finding}
L.~\color{black} Kaufman and P.~J. Rousseeuw, {\em Finding Groups in Data: An
  Introduction to Cluster Analysis}, vol.~344.
\newblock John Wiley \& Sons, 2009\color{black}.

\bibitem{hampel2011robust}
F.~R. \color{black} Hampel, E.~M. Ronchetti, P.~J. Rousseeuw, and W.~A. Stahel,
  {\em Robust {S}tatistics: {T}he {A}pproach {B}ased on {I}nfluence
  {F}unctions}, vol.~114.
\newblock John Wiley \& Sons, 2011\color{black}.

\bibitem{bai1998no}
Z.~D. Bai and J.~W. Silverstein, ``No eigenvalues outside the support of the
  limiting spectral distribution of large-dimensional sample covariance
  matrices,'' {\em Ann. Prob.}, pp.~316--345, 1998.

\bibitem{yates1995framework}
R.~D. Yates, ``A framework for uplink power control in cellular radio
  systems,'' {\em IEEE J. Sel. Areas Commun.}, vol.~13, pp.~1341--1347, Sept.
  1995.

\bibitem{horn2012matrix}
R.~A. Horn and C.~R. Johnson, {\em Matrix {A}nalysis}.
\newblock Cambridge University Press, 2012.

\end{thebibliography}

\end{document}